\documentclass{gtart_h}  

\def\ifplaintex{\expandafter\ifx\csname documentclass\endcsname\relax}

\ifplaintex 
\else
\fi

\expandafter\ifx\csname beginpicture\endcsname\relax
\expandafter\ifx\csname documentclass\endcsname\relax
\input pictex \else
\input prepictex \input pictex \input postpictex \fi\fi

\def\gt{{\mathsurround=0pt\it $\cal G\mskip-2mu$eometry \&\ 
$\cal T\!\!$opology}}        

\def\gtp{{\mathsurround=0pt\it $\cal G\mskip-2mu$eometry \&\ 
$\cal T\!\!$opology $\cal P\!$ublications}}  

\def\lognumber#1{\def\thelognumber{#1}}
\def\volumenumber#1{\def\thevolumenumber{#1}}
\def\papernumber#1{\def\thepapernumber{#1}}
\def\volumeyear#1{\def\thevolumeyear{#1}}

\def\pagenumbers#1#2{\def\startpage{#1}\def\finishpage{#2}}
\def\published#1{\def\publishdate{#1}}
\def\proposed#1{\def\theproposer{#1}}
\def\seconded#1{\def\theseconders{#1}}
\def\received#1{\def\receiveddate{#1}}

\def\accepted#1{\def\accepteddate{#1}}
\def\asciititle#1{\def\theasciititle{#1}}

\def\asciiaddress#1{\def\theasciiaddress{#1}}
\def\asciiemail#1{\def\theasciiemail{#1}}

\let\\\par\let\thelognumber\relax
\let\thevolumenumber\relax\let\thepapernumber\relax
\let\thevolumeyear\relax\let\thesamplenumber\relax\let\startpage\relax
\let\finishpage\relax\let\publishdate\relax\let\receiveddate\relax
\let\reviseddate\relax\let\accepteddate\relax\let\theasciititle\relax
\let\theasciiauthors\relax\let\theasciiaddress\relax
\let\theasciiabstract\relax
\let\theasciiemail\relax\let\theshortauthors\relax\let\theshorttitle\relax

\long\def\maketitlep{   

\count0=\startpage

\gt\hfill      
\beginpicture
\setcoordinatesystem units <0.33truein, 0.33truein> point at 2.2 0.9
\setplotsymbol ({$\cal G$})
\plotsymbolspacing=9truept
\circulararc 315 degrees from 0 1 center at 0 0
\setplotsymbol ({$\cal T$})
\circulararc 315 degrees from 1 -1 center at 1 0
\endpicture
\break
{\small\ifx\thesamplenumber\relax 
Volume \else Sample
\fi\thevolumenumber\ (\thevolumeyear)
\startpage--\finishpage\nl
Published: \publishdate}
\vglue 0.5truein plus 0.4fil minus 0.1truein

{\parskip=0pt\leftskip 0pt plus 1fil\def\\{\par\smallskip}{\ifplaintex\large
\else\Large\fi\bf\thetitle}\par\medskip}   

\vglue 0pt plus 0.1fil 

{\parskip=0pt\leftskip 0pt plus 1fil\def\\{\par}{\sc\theauthors}
\par\medskip}

\vglue 0pt plus 0.1fil 

{\small\parskip=0pt\let\newline\\
{\leftskip 0pt plus 1fil\def\\{\par}{\sl\theaddress}\par}
\expandafter\ifx\theemail\relax    
\relax\else\vglue 5pt plus 0.02fil minus 2pt\def\\{\stdspace{\rm 
and}\stdspace} 
\cl{Email:\stdspace\tt\theemail}\fi
\ifx\theurl\relax                  
\relax\else\vglue 5pt plus 0.02fil minus 2pt\def\\{\stdspace{\rm 
and}\stdspace}
\cl{URL:\stdspace\tt\theurl}\fi\par}

\vglue 7pt plus 0.3fil minus 3pt

{\bf Abstract}
\vglue 5pt plus 0.1fil minus 2pt

\theabstract

\vglue 7pt plus 0.3fil minus 3pt

{\bf AMS Classification numbers}\quad Primary:\quad \theprimaryclass

Secondary:\quad \thesecondaryclass

\vglue 5pt plus 0.3fil minus 2pt

{\bf Keywords:}\quad \thekeywords

\vglue 10pt plus 0.5fil minus 5pt

{\small  Proposed: \theproposer\hfill Received: \receiveddate\nl
Seconded: \theseconders\hfill 
\ifx\reviseddate\relax                         
Accepted: \accepteddate                        
\else
Revised: \reviseddate                          
\fi}
\eject
}       

\let\maketitlepage\maketitlep
\let\maketitle\maketitlepage

\font\phead=cmsl9 scaled 950
\font=cmsl9 scaled 1050
\font\pnum=cmbx10 scaled 913
\font\lnum=cmbx10 
\font\pfoot=cmsl9 scaled 950
\font=cmsl9 scaled 1050
\ifplaintex
\headline{\vbox to 0pt{\vskip -4.5mm\line{\small\phead\ifnum
\count0=\startpage ISSN 1364-0380 (on line)
1465-3060 (printed) \hfill {\pnum\folio}\else\ifodd\count0\def\\{ }%
\ifx\theshorttitle\relax\thetitle\else\theshorttitle\fi\hfill{\pnum\folio}
\else\def\\{ and }{\pnum\folio}\hfill\ifx\theshortauthors\relax\theauthors
\else\theshortauthors\fi\fi\fi}\vss}}
\footline{\vbox to 0pt{\vglue 0mm\line{\small\pfoot\ifnum\count0=\startpage
\copyright\ \gtp\hfill\else
\gt, Volume \thevolumenumber\ (\thevolumeyear)\hfill\fi}\vss
}}
\else
\makeatletter
\def\@oddhead{{\small\lhead\ifnum\count0=\startpage ISSN 1364-0380 (on line)
1465-3060 (printed) \hfill {\lnum\number\count0}\else\ifodd\count0
\def\\{ }\ifx\theshorttitle\relax \thetitle \else\theshorttitle\fi\hfill
{\lnum\number\count0}\else\def\\{ and }{\lnum\number\count0}
\hfill\ifx\theshortauthors\relax 
\theauthors\else\theshortauthors\fi\fi\fi}}\def\@evenhead{\@oddhead}
\def\@oddfoot{\small\lfoot\ifnum\count0=\startpage\copyright\ \gtp\hfill\else
\gt, Volume \thevolumenumber\ (\thevolumeyear)\hfill\fi}
\def\@evenfoot{\@oddfoot}
\makeatother
\fi

\newwrite\gtoutfile
\long\gdef\makeheadfile{  
{\def\\{, }\def\s{ }
\immediate\openout\gtoutfile head.xxx
\immediate\write\gtoutfile{Proxy-for: \ifx\theasciiauthors\relax
\theauthors\else\theasciiauthors\fi\s<\ifx\theasciiemail\relax\theemail\else\theasciiemail\fi>}
\immediate\write\gtoutfile{\noexpand\\}
\immediate\write\gtoutfile{Authors: \ifx\theasciiauthors\relax
\theauthors\else\theasciiauthors\fi}
{\def\\{ }\immediate\write\gtoutfile{Title: \ifx\theasciititle\relax
\thetitle\else\theasciititle\fi}}
\immediate\write\gtoutfile{Subj-class: GT or SG or MG etc}
\immediate\write\gtoutfile{MSC-class: \theprimaryclass\ifx\thesecondaryclass\relax\else, \thesecondaryclass\fi}
\immediate\write\gtoutfile{Journal-ref: Geom. Topol. \thevolumenumber
(\thevolumeyear) \startpage-\finishpage}
\immediate\write\gtoutfile{Comments: Published by Geometry and Topology at}
\immediate\write\gtoutfile{\s\s http://www.maths.warwick.ac.uk/gt/GTVol\thevolumenumber/paper\thepapernumber.abs.html}
\immediate\write\gtoutfile{\noexpand\\}
\immediate\write\gtoutfile{}
\ifx\theasciiabstract\relax
\immediate\write\gtoutfile{\theabstract}\else
\immediate\write\gtoutfile{\theasciiabstract}\fi
\immediate\write\gtoutfile{}
\immediate\write\gtoutfile{\noexpand\\}
\immediate\write\gtoutfile{}
\immediate\closeout\gtoutfile}}  

\def\maketitlepage{\maketitlep\makeheadfile}
\let\maketitle\maketitlepage

\lognumber{571}
\received{24 February 2005}
\volumenumber{9}\papernumber{36}\volumeyear{2005}
\pagenumbers{1603}{1637}   
\published{26 August 2005}
\accepted{16 August 2005}
\proposed{Yasha Eliashberg}
\seconded{Robion Kirby, Ronald Fintushel}

\usepackage{amsmath,amssymb,graphicx,psfrag}

\def\psfraga <#1,#2> #3#4{%
\psfrag {#3}{\smash{\rlap{\kern #1 \raise #2\hbox{#4}}}}}

\def\figref#1{\hyperlink{#1anchor}{Figure~\ref*{#1}}}
\def\anchor#1{\noindent\hypertarget{#1anchor}{\smash{$\phantom{99}$}}\newline}

\newtheorem{proposition}{Proposition}[section]
\newtheorem{theorem}[proposition]{Theorem}
\newtheorem{lemma}[proposition]{Lemma}
\newtheorem{corollary}[proposition]{Corollary}
\newtheorem{conjecture}[proposition]{Conjecture}

\theoremstyle{definition}
\newtheorem{definition}[proposition]{Definition}

\def\Z{\mathbb{Z}}
\def\R{\mathbb{R}}

\def\A{\mathcal{A}}

\def\I{\mathcal{I}}

\def\M{\mathcal{M}}
\def\F{\mathbb{F}}
\def\Y{\mathcal{Y}}
\def\x{X}
\def\y{Y}
\def\a{\alpha}
\def\d{\partial}

\def\Aut{\operatorname{Aut}}

\def\hom#1{\phi_{#1}}
\def\Phil#1{\Phi^L_{#1}}
\def\Phir#1{\Phi^R_{#1}}

\def\Bur#1{{\operatorname{Bur}_{#1}}}
\def\im{\operatorname{im}}
\def\ab{\operatorname{ab}}

\def\ext{\operatorname{ext}}
\def\plat{\operatorname{plat}}
\def\diag{\operatorname{diagram}}
\def\lin{\operatorname{lin}}
\def\Sym{\operatorname{Sym}}

\def\C{\mathbb{C}}
\def\CC{\mathcal{C}}
\def\H{\mathcal{H}}
\def\int{\operatorname{int}}
\def\P{\mathcal{P}}
\def\I{\mathcal{I}}

\theoremstyle{plain}

\begin{document}

\title{Knot and braid invariants from contact homology II}
\asciititle{Knot and braid invariants from contact homology II, 
with an appendix written jointly with Siddhartha Gadgil}
\author{Lenhard Ng}

\address{Department of Mathematics, Stanford
University\\Stanford, CA 94305, USA\\\medskip\\
{\rm With an appendix written jointly with
{\sc Siddhartha Gadgil}}\\\smallskip\\Stat-Math Unit, Indian
Statistical Institute\\Bangalore, India}

\asciiaddress{LN:
Department of Mathematics, Stanford
University\\Stanford, CA 94305, USA\\SG:
Stat-Math Unit, Indian
Statistical Institute\\Bangalore, India}

\asciiemail{lng@math.stanford.edu, gadgil@isibang.ac.in}
\gtemail{\mailto{lng@math.stanford.edu}, \mailto{gadgil@isibang.ac.in}}

\urladdr{http://math.stanford.edu/~lng/}

\begin{abstract}
We present a topological interpretation of knot and braid contact
homology in degree zero, in terms of cords and skein relations. This
interpretation allows us to extend the knot invariant to embedded
graphs and higher-dimensional knots. We calculate the knot invariant
for two-bridge knots and relate it to double branched covers for
general knots.

In the appendix we show that the cord ring is determined by the
fundamental group and peripheral structure of a knot and give
applications.  \end{abstract}


\maketitlepage

\section{Introduction}
\label{sec:intro}

\subsection{Main results}
\label{ssec:intro}

In \cite{I}, the author introduced invariants of knots and braid
conjugacy classes called knot and braid differential graded algebras
(DGAs). The homologies of these DGAs conjecturally give the relative
contact homology of certain natural Legendrian tori in
$5$--dimensional contact manifolds. From a computational point of
view, the easiest and most convenient way to approach the DGAs is
through the degree $0$ piece of the DGA homology, which we denoted
in \cite{I} as $HC_0$. It turns out that, unlike the full homology,
$HC_0$ is relatively easy to compute, and it gives a highly
nontrivial invariant for knots and braid conjugacy classes.

The goal of this paper is to show that $HC_0$ has a very natural
topological formulation, through which it becomes self-evident that
$HC_0$ is a topological invariant. This interpretation uses cords
and skein relations.

\begin{definition}
\label{def:cord}
Let $K\subset \R^3$ be a knot (or link). A \textit{cord} of $K$ is
any continuous path $\gamma\colon\thinspace [0,1]\rightarrow\R^3$
with $\gamma^{-1}(K) = \{0,1\}$. Denote by $\CC_K$ the set of all
cords of $K$ modulo homotopies through cords, and let $\A_K$ be the
tensor algebra over $\Z$ freely generated by $\CC_K$.
\end{definition}

\noindent In diagrams, we will distinguish between the knot and its
cords by drawing the knot more thickly than the cords.

In $\A_K$, we define \textit{skein relations} as follows:
\begin{equation}
\label{eq:skein1}
\raisebox{-0.17in}{\includegraphics[width=0.4in]{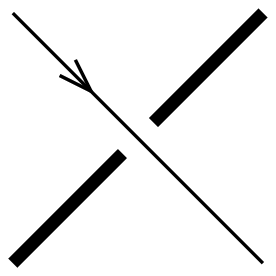}}
+
\raisebox{-0.17in}{\includegraphics[width=0.4in]{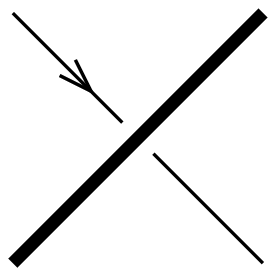}}
+
\raisebox{-0.17in}{\includegraphics[width=0.4in]{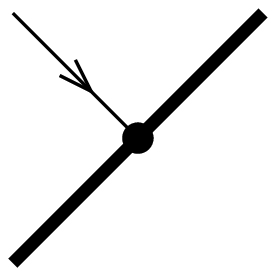}}
\cdot
\raisebox{-0.17in}{\includegraphics[width=0.4in]{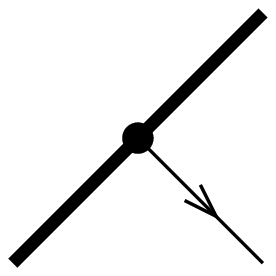}}
= 0
\end{equation}
\begin{equation}
\label{eq:skein2}
\raisebox{-0.17in}{\includegraphics[width=0.4in]{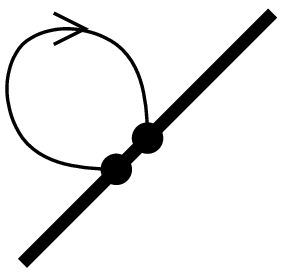}}
= -2
\end{equation}
Here, as usual, the diagrams in (\ref{eq:skein1}) are understood to
depict some local neighborhood outside of which the diagrams agree.
(The first two terms in (\ref{eq:skein1}) each show one cord, which
is split into two pieces to give the other terms.) The cord depicted
in (\ref{eq:skein2}) is any contractible cord. We write $\I_K$ as
the two-sided ideal in $\A_K$ generated by all possible skein
relations.

\begin{definition}
The \textit{cord ring} of $K$ is defined to be $\A_K/\I_K$.
\label{def:cordring}
\end{definition}

\noindent It is clear that the cord ring yields a topological
invariant of the knot; for a purely homotopical definition of the
cord ring, in joint work with S Gadgil, see the Appendix. However,
it is not immediately obvious that this ring is small enough to be
manageable (for instance, finitely generated), or large enough to be
interesting. The main result of this paper is the following.

\begin{theorem}
\label{thm:mainknot}
The cord ring of $K$ is isomorphic to the degree $0$ knot contact
homology $HC_0(K)$.
\end{theorem}

\noindent For the definition of $HC_0$, see
Section~\ref{ssec:background}.

\begin{figure}[ht!]\anchor{fig:trefoil}
\centerline{
\small
\psfrag {g1}{$\gamma_1$}
\psfrag {g2}{$\gamma_2$}
\psfrag {g3}{$\gamma_3$}
\psfrag {g4}{$\gamma_4$}
\psfraga <-2pt,0pt> {g5}{$\gamma_5$}
\includegraphics[width=2.5in]{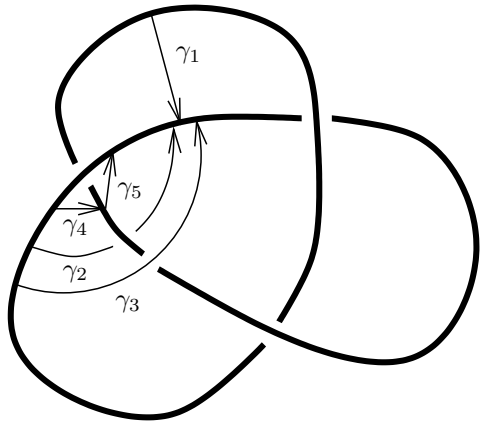}
}
\caption{
A trefoil, with a number of its cords
}
\label{fig:trefoil}
\end{figure}

As an example, consider the trefoil $3_1$ in
\figref{fig:trefoil}. By keeping the ending point fixed and
swinging the beginning point around the trefoil, we see that
$\gamma_1$ is homotopic to both $\gamma_2$ and $\gamma_5$;
similarly, $\gamma_4$ is homotopic to $\gamma_1$ (move both
endpoints counterclockwise around the trefoil), and $\gamma_3$ is
homotopic to a trivial loop. On the other hand, skein relation
(\ref{eq:skein1}) implies that, in $HC_0(3_1) = \A_K/\I_K$, we have
$\gamma_2+\gamma_3+\gamma_4\gamma_5=0$, while skein relation
(\ref{eq:skein2}) gives $\gamma_3=-2$. We conclude that
\[
0 = \gamma_4\gamma_5 + \gamma_2 + \gamma_3 = \gamma_1^2 + \gamma_1
-2.
\]
In fact, it turns out that $HC_0(3_1)$ is generated by $\gamma_1$
with relation $\gamma_1^2+\gamma_1-2$; see Section~\ref{ssec:plat}.

We can extend our definitions to knots in arbitrary $3$--manifolds.
In particular, a braid $B$ in the braid group $B_n$ yields a knot in
the solid torus $D^2\times S^1$, and the isotopy class of this knot
depends only on the conjugacy class of $B$. If we define $\A_B$ and
$\I_B$ as above, with $B$ as the knot in $D^2\times S^1$, then we
have the following analogue of Theorem~\ref{thm:mainknot}.

\begin{theorem}
\label{thm:mainbraid}
The cord ring $\A_B/\I_B$ is isomorphic to the degree $0$ braid
contact homology $HC_0(B)$.
\end{theorem}

There is also a version of the cord ring involving unoriented cords.
The \textit{abelian cord ring} for a knot is the commutative ring
generated by unoriented cords, modulo the skein relations
(\ref{eq:skein1}) and (\ref{eq:skein2}). In other words, it is the
abelianization of the cord ring modulo identifying cords and their
orientation reverses. Analogues of Theorems~\ref{thm:mainknot}
and~\ref{thm:mainbraid} then state that the abelian cord rings of a
knot $K$ or a braid $B$ are isomorphic to the rings $HC_0^{\ab}(K)$
or $HC_0^{\ab}(B)$ (see Section~\ref{ssec:background}).

The cord ring formulation of $HC_0$ is useful in several ways
besides its intrinsic interest. In \cite{I}, we demonstrated how to
calculate $HC_0$ for a knot, via a closed braid presentation of the
knot. Using the cord ring, we will see how to calculate $HC_0$
instead in terms of either a plat presentation or a knot diagram,
which is more efficient in many examples. In particular, we can
calculate $HC_0$ for all $2$--bridge knots
(Theorem~\ref{thm:2bridge}). The cord ring can also be applied to
find lower bounds for the number of minimal-length chords of a knot.

It was demonstrated in \cite{I} that $HC_0$ is related to the
determinant of the knot. An intriguing application of the cord
formalism is a close connection between the abelian cord ring
$HC_0^{\ab}$ and the $SL_2(\C)$ character variety of the double
branched cover of the knot (Proposition~\ref{prop:char}).

In addition, the cord ring is defined in much more generality than
just for knots and braids. We have already mentioned that it gives a
topological invariant of knots in any $3$--manifold. It also extends
to embedded graphs in $3$--manifolds, for which it gives an invariant
under neighborhood equivalence, and to knots in higher dimensions.

We now outline the paper. Section~\ref{ssec:background} is included
for completeness, and contains the definitions of knot and braid
contact homology. In Section~\ref{sec:braidrep}, we examine the
braid representation used to define contact homology. This
representation was first introduced by Magnus in relation to
automorphisms of free groups; our geometric interpretation, which is
reminiscent of the ``forks'' used by Krammer \cite{Kra} and Bigelow
\cite{Big} to prove linearity of the braid groups, is crucial to the
identification of the cord ring with $HC_0$. We extend this
geometric viewpoint in Section~\ref{sec:cords} and use it to prove
Theorems~\ref{thm:mainknot} and~\ref{thm:mainbraid}. In
Section~\ref{sec:plat}, we discuss how to calculate the cord ring in
terms of either plats or knot diagrams, with a particularly simple
answer for $4$--plats. Section~\ref{sec:remarks} discusses some
geometric consequences, including connections to double branched
covers and an extension of the cord ring to the graph invariant
mentioned previously. The Appendix, written with S Gadgil, gives a
group-theoretic formulation for the cord ring, and discusses an
extension of the cord ring to a nontrivial invariant of codimension
$2$ submanifolds in any manifold.

\subsection*{Acknowledgements}
 I am grateful to Dror Bar-Natan, Tobias Ekholm, Yasha
Eliashberg, Siddhartha Gadgil, and Justin Roberts for interesting
and useful conversations, and to Stanford University and the
American Institute of Mathematics for their hospitality. The work
for the Appendix was done at the June 2003 workshop on holomorphic
curves and contact geometry in Berder, France. This work is
supported by a Five-Year Fellowship from the American Institute of
Mathematics.

\subsection{Background material}
\label{ssec:background}

We recall the definitions of degree $0$ braid and knot contact
homology from \cite{I}. Let $\A_n$ denote the tensor algebra over
$\Z$ generated by $n(n-1)$ generators $a_{ij}$ with $1\leq i,j\leq
n$, $i\neq j$. There is a representation $\phi$ of the braid group
$B_n$ as a group of algebra automorphisms of $\A_n$, defined on
generators $\sigma_k$ of $B_n$ by:
\[
\hom{\sigma_k}\colon\thinspace \left\{
\begin{array}{ccll}
a_{ki} & \mapsto & -a_{k+1,i} - a_{k+1,k}a_{ki} & i\neq k,k+1 \\
a_{ik} & \mapsto & -a_{i,k+1} - a_{ik}a_{k,k+1} & i\neq k,k+1 \\
a_{k+1,i} & \mapsto & a_{ki} & i\neq k,k+1 \\
a_{i,k+1} & \mapsto & a_{ik} & i \neq k,k+1 \\
a_{k,k+1} & \mapsto & a_{k+1,k} & \\
a_{k+1,k} & \mapsto & a_{k,k+1} & \\
a_{ij} & \mapsto & a_{ij} & i,j \neq k,k+1
\end{array}
\right.
\]
In general, we denote the image of $B\in B_n$ in $\Aut\A_n$ by
$\hom{B}$.

\begin{definition}
For $B\in B_n$, the degree $0$ braid contact homology is defined by
$HC_0(B) = \A_n/\im(1-\hom{B})$, where $\im(1-\hom{B})$ is the
two-sided ideal in $\A_n$ generated by the image of the map
$1-\hom{B}$.
\end{definition}

To define knot contact homology, we need a bit more notation.
Consider the map $\phi^{\ext}$ given by the composition $B_n
\hookrightarrow B_{n+1} \stackrel{\phi}{\rightarrow} \Aut\A_{n+1}$,
where the inclusion simply adds a trivial strand labeled $*$ to any
braid. Since $*$ does not cross the other strands, we can express
$\phi^{\ext}_B(a_{i*})$ as a linear combination of $a_{j*}$ with
coefficients in $\A_n$, and similarly for $\phi^{\ext}_B(a_{j*})$.
More concretely, for $B\in B_n$, define matrices $\Phil{B},\Phir{B}$
by
\[
\hom{B}^{\ext}(a_{i*}) = \sum_{j=1}^n (\Phil{B})_{ij} a_{j*}
\hspace{0.25in} \textrm{and} \hspace{0.25in}
 \hom{B}^{\ext}(a_{*j}) =
\sum_{i=1}^n a_{*i} (\Phir{B})_{ij}.
\]
Also, define for convenience the matrix $A = (a_{ij})$; here and
throughout the paper, we set $a_{ii}=-2$ for any $i$.

\begin{definition}
If $K$ is a knot in $\R^3$, let $B\in B_n$ be a braid whose closure
is $K$. Then the degree $0$ knot contact homology of $K$ is defined
by $HC_0(K) = \A_n/I$, where $I$ is the two-sided ideal in $\A_n$
generated by the entries of the matrices $A-\Phil{B}\cdot A$ and
$A-A\cdot\Phir{B}$. Up to isomorphism, this depends only on $K$ and
not on the choice of $B$.
\label{def:knotcordring}
\end{definition}

Finally, the abelian versions of $HC_0$ are defined as follows:
$HC_0^{\ab}(B)$ and $HC_0^{\ab}(K)$ are the abelianizations of
$HC_0(B)$ and $HC_0(K)$, modulo setting $a_{ij}=a_{ji}$ for all
$i,j$.

The main results of \cite{I} state, in part, that $HC_0(B)$ and
$HC_0^{\ab}(B)$ are invariants of the conjugacy class of $B$, while
$HC_0(K)$ and $HC_0^{\ab}(K)$ are knot invariants. As mentioned in
Section~\ref{ssec:intro}, these results follow directly from
Theorems~\ref{thm:mainknot} and~\ref{thm:mainbraid} here.

\section{Braid representation revisited}
\label{sec:braidrep}

The braid representation $\phi$ was introduced and studied, in a
slightly different form, by Magnus \cite{Mag} and then Humphries
\cite{Hum}, both of whom treated it essentially algebraically. In
this section, we will give a geometric interpretation for $\phi$.
Our starting point is the well-known expression of $B_n$ as a
mapping class group.

Let $D$ denote the unit disk in $\C$, and let $P=\{p_1,\ldots,p_n\}$
be a set of distinct points (``punctures'') in the interior of $D$.
We will choose $P$ such that $p_i\in\R$ for all $i$, and
$p_1<p_2<\cdots<p_n$; in figures, we will normally omit drawing the
boundary of $D$, and we depict the punctures $p_i$ as dots. Write
$\H(D,P)$ for the set of orientation-preserving homeomorphisms $h$
of $D$ satisfying $h(P)=P$ and $h|_{\d D} = \operatorname{id}$, and
let $\H_0(D,P)$ be the identity component of $\H(D,P)$. Then $B_n =
\H(D,P)/\H^0(D,P)$, the mapping class group of $(D,P)$ (for
reference, see \cite{Bir}). We will adopt the convention that the
generator $\sigma_k\in B_n$ interchanges the punctures $p_k,p_{k+1}$
in a counterclockwise fashion while leaving the other punctures
fixed.

\begin{figure}[ht!]\anchor{fig:arcs}
\centerline{\small
\psfrag {a}{$a_{ij}$}
\psfrag {pi}{$p_{i}$}
\psfrag {pj}{$p_{j}$}
\includegraphics[width=3.5in]{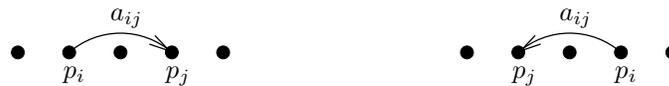}
}
\caption{
The arcs $\gamma_{ij}$ for $i<j$ (left) and $i>j$ (right)
}
\label{fig:arcs}
\end{figure}

\begin{definition}
An (oriented) \textit{arc} is an embedding $\gamma\colon\thinspace
[0,1]\rightarrow\int(D)$ such that $\gamma^{-1}(P)=\{0,1\}$. We
denote the set of arcs modulo isotopy by $\P_n$. For $1\leq i,j\leq
n$ with $i\neq j$, we define $\gamma_{ij}\in\P_n$ to be the arc from
$p_i$ to $p_j$ which remains in the upper half plane; see
\figref{fig:arcs}.
\end{definition}

\noindent The terminology derives from \cite{Kra}, where
(unoriented) arcs are used to define the Lawrence--Krammer
representation of $B_n$. Indeed, arcs are central to the proofs by
Bigelow and Krammer that this representation is faithful. We remark
that it might be possible to recover Lawrence--Krammer from the
algebra representation $\phi^{\ext}$, using arcs as motivation.

The braid group $B_n$ acts on $\P_n$ via the identification with the
mapping class group. The idea underlying this section is that there
is a map from $\P_n$ to $\A_n$ under which this action corresponds
to the representation $\phi$.

\begin{proposition}
\label{prop:psi}
There is an unique map $\psi\colon\thinspace \P_n\rightarrow\A_n$
satisfying the following properties:
\begin{enumerate}
\item {\rm Equivariance}\qua $\psi(B\cdot \gamma) = \hom{B}(\psi(\gamma))$ for
any $B\in B_n$ and $\gamma\in\P_n$, where $B\cdot\gamma$ denotes the
action of $B$ on $\gamma$.
\item {\rm Normalization}\qua $\psi(\gamma_{ij})=a_{ij}$ for all $i,j$.
\end{enumerate}
\end{proposition}

\begin{proof}
Since the action of $B_n$ on $\P_n$ is transitive, we define
$\psi(\gamma)$ by choosing any $B_\gamma\in B_n$ for which
$B_\gamma\gamma_{12} = \gamma$ and then setting $\psi(\gamma) =
\hom{B_{\gamma}}(a_{12})$. (This shows that $\psi$, if it exists,
must be unique.) First assume that this yields a well-defined map.
Then for $B\in B_n$, we have $B\cdot\gamma =
BB_\gamma\cdot\gamma_{12}$, and so
\[
\psi(B\cdot\gamma) = \hom{BB_{\gamma}}(a_{12}) =
\hom{B}(\hom{B_{\gamma}}(a_{12})) = \hom{B}\psi(\gamma).
\]
In addition, if $i<j$, then $(\sigma_{i-1}^{-1}\cdots
\sigma_1^{-1})(\sigma_{j-1}^{-1}\cdots \sigma_2^{-1})$ maps
$\gamma_{12}$ to $\gamma_{ij}$, while $\hom{\sigma_{i-1}^{-1}\cdots
\sigma_1^{-1}\sigma_{j-1}^{-1}\cdots \sigma_2^{-1}}(a_{12}) =
a_{ij}$; if $i>j$, then $(\sigma_{j-1}^{-1}\cdots
\sigma_1^{-1})(\sigma_{i-1}^{-1}\cdots \sigma_2^{-1})\sigma_1$ maps
$\gamma_{12}$ to $\gamma_{ij}$, while $\hom{\sigma_{j-1}^{-1}\cdots
\sigma_1^{-1}\sigma_{i-1}^{-1}\cdots \sigma_2^{-1}\sigma_1}(a_{12})
= a_{ij}$.

We now only need to show that the above definition of $\psi$ is
well-defined. By transitivity, it suffices to show that if
$B\cdot\gamma_{12}=\gamma_{12}$, then $\hom{B}(a_{12})=a_{12}$. Now
if $B\cdot\gamma_{12} = \gamma_{12}$, then $B$ preserves a
neighborhood of $\gamma_{12}$; if we imagine contracting this
neighborhood to a point, then $B$ becomes a braid in $B_{n-1}$ which
preserves this new puncture. Now the subgroup of $B_{n-1}$ which
preserves the first puncture (ie, whose projection to the
symmetric group $S_{n-1}$ keeps $1$ fixed) is generated by
$\sigma_k$, $2\leq k \leq n-2$, and
$(\sigma_1\sigma_2\cdots\sigma_{k-1})(\sigma_{k-1}\cdots\sigma_2\sigma_1)$,
$2\leq k \leq n-1$. It follows that the subgroup of braids $B\in
B_n$ which preserve $\gamma_{12}$ is generated by $\sigma_1^2$
(which revolves $\gamma_{12}$ around itself); $\sigma_k$ for $3 \leq
k \leq n-1$; and
\[
\tau_k = (\sigma_2\sigma_1)(\sigma_3\sigma_2)\cdots
(\sigma_{k-1}\sigma_{k-2})(\sigma_{k-2}\sigma_{k-1})\cdots
(\sigma_2\sigma_3)(\sigma_1\sigma_2)
\]
for $3 \leq k \leq n$. But $\phi_{\sigma_1^2}$ and $\phi_{\sigma_k}$
clearly preserve $a_{12}$, while $\phi_{\tau_k}$ preserves $a_{12}$
because $\hom{\sigma_i\sigma_{i+1}}(a_{i,i+1}) = a_{i+1,i+2}$ and
$\hom{\sigma_{i+1}\sigma_i}(a_{i+1,i+2}) = a_{i,i+1}$ for $1 \leq i
\leq n-2$.
\end{proof}

The map $\psi$ satisfies a skein relation analogous to the skein
relation from Section~\ref{sec:intro}.

\begin{proposition}
\label{prop:psiskein}
The following skein relation holds for arcs:
\begin{equation}
\label{eq:psiskein}
\psi(\raisebox{-0.15in}{\includegraphics[width=0.6in]{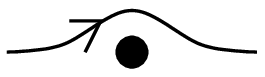}})
+
\psi(\raisebox{-0.15in}{\includegraphics[width=0.6in]{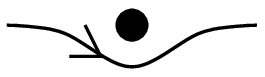}})
+
\psi(\raisebox{-0.15in}{\includegraphics[width=0.6in]{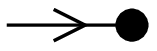}})
\psi(\raisebox{-0.15in}{\includegraphics[width=0.6in]{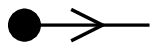}})
= 0.
\end{equation}
\end{proposition}

\begin{proof}
By considering the concatenation of the two arcs involved in the
product in the above identity, which are disjoint except for one
shared endpoint, we see that there is some element of $B_n$ which
maps the two arcs to $\gamma_{12}$ and $\gamma_{23}$. Since $\psi$
is $B_n$--equivariant, it thus suffices to establish the identity
when the two arcs are $\gamma_{12}$ and $\gamma_{23}$. In this case,
the other two arcs in the identity are $\gamma_{13}$ and $\gamma$,
where $\gamma$ is a path joining $p_1$ to $p_3$ lying in the lower
half plane. But then $\gamma = \sigma_2 \cdot \gamma_{12}$, and
hence by normalization and equivariance,
\[
\psi(\gamma_{13}) + \psi(\gamma) + \psi(\gamma_{12})
\psi(\gamma_{23}) = a_{13} + \hom{\sigma_2}(a_{12}) + a_{12} a_{23}
= 0,
\]
as desired.
\end{proof}

Rather than defining $\psi$ in terms of $\phi$, we could imagine
first defining $\psi$ via the normalization of
Proposition~\ref{prop:psi} and the skein relation
(\ref{eq:psiskein}), and then defining $\phi$ by $\hom{B}(a_{ij}) =
\psi(B\cdot \gamma_{ij})$. For instance, (\ref{eq:psiskein}) implies
that
\begin{eqnarray*}
\psi(\sigma_1\cdot\gamma_{13}) &=&
\psi(\raisebox{-0.11in}{\footnotesize
\psfrag {p1}{$p_1$}
\psfrag {p2}{$p_2$}
\psfrag {p3}{$p_3$}
\includegraphics[height=0.35in]{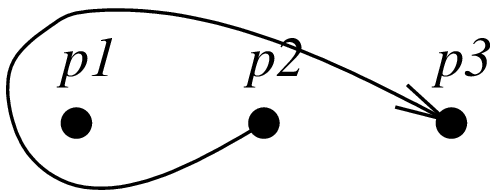}}) \\
&=& -
\psi(\raisebox{-0.16in}{\footnotesize
\psfrag {p1}{$p_1$}
\psfrag {p2}{$p_2$}
\psfrag {p3}{$p_3$}
\includegraphics[height=0.45in]{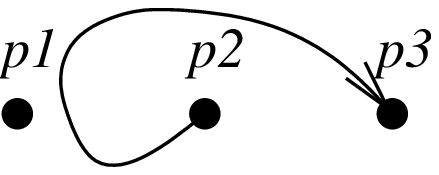}})
-\psi(\raisebox{-0.16in}{\footnotesize
\psfrag {p1}{$p_1$}
\psfrag {p2}{$p_2$}
\psfrag {p3}{$p_3$}
\includegraphics[height=0.45in]{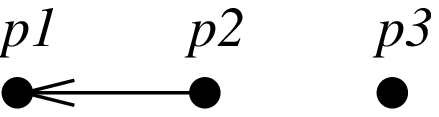}})
\psi(\raisebox{-0.16in}{\footnotesize
\psfrag {p1}{$p_1$}
\psfrag {p2}{$p_2$}
\psfrag {p3}{$p_3$}
\includegraphics[height=0.45in]{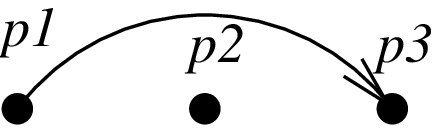}}) \\
&=& -a_{23}-a_{21}a_{13},
\end{eqnarray*}
which gives $\hom{\sigma_1}(\gamma_{13}) = -a_{23}-a_{21}a_{13}$.

\begin{proposition}
\label{prop:suffice}
The skein relation of Proposition~\ref{prop:psiskein} and the
normalization property of Proposition~\ref{prop:psi} suffice to
define the map $\psi\colon\thinspace \P_n\rightarrow\A_n$.
\end{proposition}

Before proving Proposition~\ref{prop:suffice}, we need to introduce
some notation.
\begin{definition}
An arc $\gamma\in\P_n$ is in \textit{standard form} if its image in
$D$ consists of a union of semicircles centered on the real line,
each contained in either the upper half plane or the lower half
plane. An arc is in \textit{minimal standard form} if it is in
standard form, and either it lies completely in the upper half
plane, or each semicircle either contains another semicircle in the
same half plane nested inside of it, or has a puncture along its
diameter (not including endpoints).
\end{definition}

\begin{figure}[ht!]\anchor{fig:standard}
\centerline{
\includegraphics[width=3.5in]{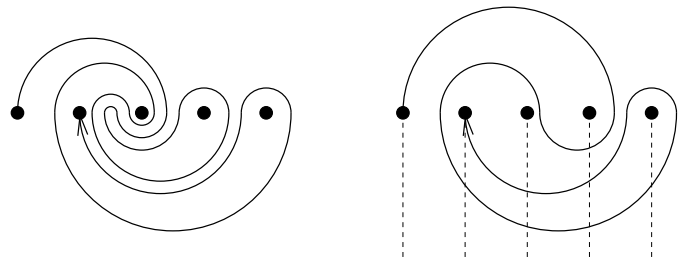}
}
\caption{
An arc in standard form (left), and the corresponding minimal
standard form (right). The dashed lines are used to calculate the
height of the arc in minimal standard form, which is $8$ in this
case.
}
\label{fig:standard}
\end{figure}

See \figref{fig:standard} for examples. It is easy to see that
any arc can be perturbed into standard form while fixing all
intersections with the real line, and any arc in standard form can
be isotoped to an arc in minimal standard form.

Define the height $h$ of any arc as follows: for each puncture, draw
a ray starting at the puncture in the negative imaginary direction,
and count the number of (unsigned) intersections of this ray with
the arc, where an endpoint of the arc counts as half of a point; the
height is the sum of these intersection numbers over all punctures.
(See \figref{fig:standard}. Strictly speaking, $h$ is only
defined for arcs which are not tangent to the rays anywhere outside
of their endpoints, but this will not matter.) An isotopy sending
any arc to an arc in minimal standard form does not increase height;
that is, minimal standard form minimizes height for any isotopy
class of arcs.

The following is the key result which allows us to prove
Proposition~\ref{prop:suffice}, as well as faithfulness results for
$\phi$.

\begin{lemma}
\label{lemma:induct}
Let $\gamma$ be a minimal standard arc with $h(\gamma)>1$. Then
there are minimal standard arcs $\gamma',\gamma_1,\gamma_2$ with
$h(\gamma')<h(\gamma) = h(\gamma_1)+h(\gamma_2)$ related by the
skein relation $\psi(\gamma) =
-\psi(\gamma')-\psi(\gamma_1)\psi(\gamma_2)$.
\end{lemma}

\begin{proof}
Define a \textit{turn} of $\gamma$ to be any point on $\gamma$
besides the endpoints for which the tangent line to $\gamma$ is
vertical (parallel to the imaginary axis); note that all turns lie
on the real line. We consider two cases.

If $\gamma$ has $0$ turns or $1$ turn, then by minimality, it
contains a semicircle in the lower half plane whose diameter
includes a puncture distinct from the endpoints of $\gamma$. We can
use the skein relation to push the semicircle through this puncture.
When $\gamma$ is pushed to pass through the puncture, it splits into
two arcs $\tilde{\gamma}_1,\tilde{\gamma}_2$ whose heights sum to
$h(\gamma)$; after it passes the puncture, it gives an arc
$\tilde{\gamma}'$ whose height is $h(\gamma)-1$. When we isotop all
of these arcs to minimal standard forms $\gamma_1,\gamma_2,\gamma'$,
the height of $\tilde{\gamma}'$ does not increase, while the heights
of $\tilde{\gamma}_1,\tilde{\gamma}_2$ are unchanged. The lemma
follows in this case.

Now suppose that $\gamma$ has at least $2$ turns. Let $q$ be a turn
representing a local maximum of the real part of $\gamma$, let $p$
be the closest puncture to the left of $q$ (ie, the puncture whose
value in $\R$ is greatest over all punctures less than $q$); by
replacing $q$ if necessary, we can assume that $q$ is the closest
turn to the right of $p$. Now there are two semicircles in $\gamma$
with endpoint at $q$; by minimality, the other endpoints of these
semicircles are to the left of $p$. We can thus push $\gamma$
through $p$ so that the turn $q$ passes across $p$, and argue as in
the previous case, unless $p$ is an endpoint of $\gamma$.

Since we can perform a similar argument for a turn representing a
local minimum, we are done unless the closest puncture to the
left/right of any max/min turn (respectively) is an endpoint of
$\gamma$. We claim that this is impossible. Label the endpoints of
$\gamma$ as $p_1<p_2$, and traverse $\gamma$ from $p_1$ to $p_2$. It
is easy to see from minimality that the first turn we encounter must
be to the right of $p_2$, while the second must be to the left of
$p_1$. This forces the existence of a third turn to the right of
$p_2$, and a fourth to the left of $p_1$, and so forth, spiraling
out indefinitely and making it impossible to reach $p_2$.
\end{proof}

\begin{proof}[Proof of Proposition~\ref{prop:suffice}]
By Lemma~\ref{lemma:induct}, we can use the skein relation to
express (the image under $\psi$ of) any minimal standard arc of
height at least $2$ in terms of minimal standard arcs of strictly
smaller height, since any arc has height at least $1$. The
normalization condition defines the image under $\psi$ of arcs of
height $1$, and the proposition follows.
\end{proof}

We now examine the question of the faithfulness of $\phi$. Define
the degree operator on $\A_n$ as usual: if $v\in\A_n$, then $\deg v$
is the largest $m$ such that there is a monomial in $v$ of the form
$ka_{i_1j_1}a_{i_2j_2}\cdots a_{i_mj_m}$.

\begin{proposition}
\label{prop:deg}
For $\gamma\in\P_n$ a minimal standard arc, $\deg\psi(\gamma) =
h(\gamma)$.
\end{proposition}

\begin{proof}
This is an easy induction on the height of $\gamma$, using
Lemma~\ref{lemma:induct}. If $h(\gamma)=1$, then
$\gamma=\gamma_{ij}$ for some $i,j$, and so $\psi(\gamma)=a_{ij}$
has degree $1$. Now assume that the assertion holds for
$h(\gamma)\leq m$, and consider $\gamma$ with $h(\gamma)=m+1$. With
notation as in Lemma~\ref{lemma:induct}, we have
$h(\gamma'),h(\gamma_1),h(\gamma_2) \leq m$, and so
$\deg(\psi(\gamma_1)\psi(\gamma_2)) = h(\gamma_1)+h(\gamma_2)=m+1$
while $\deg(\psi(\gamma')) \leq m$. It follows that
$\deg\psi(\gamma)=m+1$, as desired.
\end{proof}

\begin{corollary}
The map $\psi\colon\thinspace \P_n\rightarrow\A_n$ is injective.
\label{cor:injective}
\end{corollary}

\begin{proof}
Suppose $\gamma,\gamma'\in\P_n$ satisfy
$\psi(\gamma)=\psi(\gamma')$. Since $\psi$ is $B_n$--equivariant and
$B_n$ acts transitively on $\P_n$, we may assume that
$\gamma'=\gamma_{12}$. We may further assume that $\gamma$ is a
minimal standard arc; then by Proposition~\ref{prop:deg},
$h(\gamma)=h(\gamma_{12})=1$, and so $\gamma$ is isotopic to
$\gamma_{ij}$ for some $i,j$. Since
$\psi(\gamma_{12})=\psi(\gamma_{ij}) = a_{ij}$, we conclude that
$i=1,j=2$, and hence $\gamma$ is isotopic to $\gamma_{12}$.
\end{proof}

We next address the issue of faithfulness. Recall from \cite{Hum} or
by direct computation that $\phi\colon\thinspace
B_n\rightarrow\Aut(\A_n)$ is not a faithful representation; its
kernel has been shown in \cite{Hum} to be the center of $B_n$, which
is generated by $(\sigma_1\cdots\sigma_{n-1})^n$. However, the
extension $\phi^{\ext}$ discussed in Section~\ref{ssec:background}
is faithful, as was first shown in \cite{Mag}.

To interpret $\phi^{\ext}$ in the mapping class group picture, we
introduce a new puncture $*$, which we can think of as lying on the
boundary of the disk, and add this to the usual $n$ punctures; $B_n$
now acts on this punctured disk in the usual way, in particular
fixing $*$. The generators of $\A_{n+1}$ not in $\A_n$ are of the
form $a_{i*},a_{*i}$, with corresponding arcs
$\gamma_{i*},\gamma_{*i}\subset D$. Although we have previously
adopted the convention that all punctures lie on the real line, we
place $*$ at the point $\sqrt{-1}\in D$ for convenience, with
$\gamma_{i*},\gamma_{*i}$ the straight line segments between $*$ and
puncture $p_i\in\R$. As in Propositions~\ref{prop:psi} and
~\ref{prop:psiskein}, there is a map $\psi^{\ext}\colon\thinspace
\P_{n+1}\rightarrow\A_{n+1}$ defined by the usual skein relation
(\ref{eq:psiskein}), or alternatively by $\psi^{\ext}(B\cdot\gamma)
= \phi^{\ext}_B(\psi(\gamma))$ for any $B\in B_n$ and
$\gamma\in\P_{n+1}$.

We are now in a position to give a geometric proof of the
faithfulness results from \cite{Hum} and \cite{Mag}.

\begin{proposition}\label{prop:faithful}
{\rm\cite{Hum,Mag}}\qua
The map $\phi^{\ext}$ is faithful, while the kernel of $\phi$ is the
center of $B_n$, $\{(\sigma_1\cdots\sigma_{n-1})^{nm}\,|\,m\in\Z\}$.
\end{proposition}

\begin{proof}
We first show that $\phi^{\ext}$ is faithful. Suppose that $B\in
B_n$ satisfies $\phi^{\ext}_B = 1$. Then, in particular,
$\psi^{\ext}(B\cdot a_{*i}) = \phi^{\ext}_B(a_{*i}) = a_{*i}$, and
so by Corollary~\ref{cor:injective}, the homeomorphism $f_B$ of $D$
determined by $B$ sends $\gamma_{*i}$ to an arc isotopic to
$\gamma_{*i}$ for all $i$. This information completely determines
$f_B$ up to isotopy and implies that $B$ must be the identity braid
in $B_n$. (One can imagine cutting open the disk along the arcs
$\gamma_{*i}$ to obtain a puncture-free disk on which $f_B$ is the
identity on the boundary; it follows that $f_B$ must be isotopic to
the identity map.) See \figref{fig:faithful}.

\begin{figure}[ht!]\anchor{fig:faithful}
\centerline{\small
\psfrag {p1}{$p_1$}
\psfrag {p2}{$p_2$}
\psfrag {pn}{$p_n$}
\includegraphics[width=3.5in]{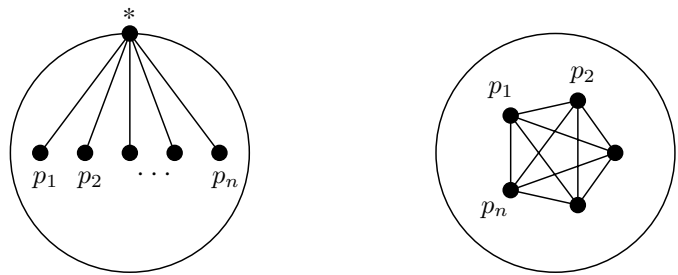}
}
\caption{
Proof of Proposition~\ref{prop:faithful}. If $\phi^{\ext}_B = 1$,
then $B$ preserves all of the arcs in the left diagram, oriented in
either direction; if $\hom{B} = 1$, then $B$ preserves all of the
arcs in the right diagram.
}
\label{fig:faithful}
\end{figure}

A similar argument can be used for computing $\ker \phi$. Rearrange
the punctures $p_1,\ldots,p_n$ in a circle, so that $\gamma_{ij}$
becomes the line segment from $p_i$ to $p_j$ for all $i,j$ (see
\figref{fig:faithful}). If $B\in\ker\phi$, then the
homeomorphism $f_B$ determined by $B$ sends each $\gamma_{ij}$ to an
arc isotopic to $\gamma_{ij}$. We may assume without loss of
generality that $f_B$ actually preserves each $\gamma_{ij}$; then,
by deleting the disk bounded by
$\gamma_{12},\gamma_{23},\ldots,\gamma_{n-1,n},\gamma_{n1}$, we can
view $f_B$ as a homeomorphism of the annulus which is the identity
on both boundary components. For any such homeomorphism, there is an
$m\in\Z$ such that the homeomorphism is isotopic to the map which
keeps the outside boundary fixed and rotates the rest of the annulus
progressively so that the inside boundary is rotated by $m$ full
revolutions. This latter map corresponds to the braid
$(\sigma_1\cdots\sigma_{n-1})^{nm}$; the result follows.
\end{proof}

For future use, we can also give a geometric proof of a result from
\cite{I}.

\begin{proposition}{\rm \cite[Proposition 4.7]{I}}
\label{prop:matrix}
We have the matrix identity $(\hom{B}(a_{ij}))$ $ = A - \Phil{B} \cdot
A \cdot \Phir{B}$.
\end{proposition}

\begin{proof}
We can write $\gamma_{ij}$ as the union of the arcs $\gamma_{i*}$
and $\gamma_{*j}$, which are disjoint except at $*$. Thus
$B\cdot\gamma_{ij}$ is the union of the arcs $B\cdot\gamma_{i*}$ and
$B\cdot\gamma_{*j}$. Now, by the definition of $\Phil{B}$, we can
write
\[
\psi(B\cdot\gamma_{i*}) = \phi_B^{\ext}(a_{i*}) = \sum_k
(\Phil{B})_{ik} a_{k*}
\]
and similarly $\psi(B\cdot\gamma_{*j}) = \sum_l a_{*l}
(\Phir{B})_{lj}$. Since the union of the arcs $\gamma_{k*}$ and
$\gamma_{*l}$ is $\gamma_{kl}$, it follows that
\[
\psi(B\cdot\gamma_{ij}) = \sum_{k,l} (\Phil{B})_{ik} a_{kl}
(\Phir{B})_{lj}.
\]
Assembling this identity in matrix form gives the proposition.
\end{proof}

\section{Cords and the cord ring}
\label{sec:cords}

\subsection{Cords in $(D,P)$}
\label{ssec:cordmap}
It turns out that the map $\psi$ on (embedded) arcs can be extended
to paths which are merely immersed. This yields another description
of $\psi$, independent from the representation $\phi$. We give this
description in this section, and use it in
Section~\ref{ssec:thmproofs} to prove Theorems~\ref{thm:mainknot}
and~\ref{thm:mainbraid}.

\begin{definition}
\label{def:cordplane}
A \textit{cord} in $(D,P)$ is a continuous map
$\gamma\colon\thinspace [0,1]\rightarrow\int(D)$ with
$\gamma^{-1}(P)=\{0,1\}$. (In particular, $\gamma(0)$ and
$\gamma(1)$ are not necessarily distinct.) We denote the set of
cords in $(D,P)$, modulo homotopy through cords, by $\tilde{\P}_n$.
\end{definition}

Given a cord $\gamma$ in $(D,P)$ with $\gamma(0)=p_i$ and
$\gamma(1)=p_j$, there is a natural way to associate an element
$\x(\gamma)$ of $\F_n$, the free group on $n$ generators
$x_1,\ldots,x_n$, which we identify with $\pi_1(D\setminus P)$ by
setting $x_m$ to be the counterclockwise loop around $p_m$.
Concatenate $\gamma$ with the arc $\gamma_{ji}$; this gives a loop,
for which we choose a base point on $\gamma_{ji}$. (If $i=j$, then
$\gamma$ already forms a loop, and we can choose any base point on
$\gamma$ in a neighborhood of $p_i=p_j$.) If we push this loop off
of the points $p_i$ and $p_j$, we obtain a based loop $\x(\gamma)
\in \pi_1(D\setminus P) = \F_n$. It is important to note that
$\x(\gamma)$ is only well-defined up to multiplication on the left
by powers of $x_i$, and on the right by powers of $x_j$.

We wish to extend the map $\psi$ to $\tilde{\P}_n$. To do this, we
introduce an auxiliary tensor algebra $\Y_n$ over $\Z$ on $n$
generators $y_1,\ldots,y_n$. There is a map $\y\colon\thinspace
\F_n\rightarrow\Y_n/\langle y_1^2+2y_1,\ldots,y_n^2+2y_n\rangle$
defined on generators by $\y(x_i) = \y(x_i^{-1}) = -1-y_i$, and
extended to $\F_n$ in the obvious way: $\y(x_{i_1}^{k_1} \cdots
x_{i_m}^{k_m}) = (-1-y_{i_1})^{k_1} \cdots (-1-y_{i_m})^{k_m}$. This
is well-defined since $\y(x_i)\y(x_i^{-1}) = \y(x_i^{-1})\y(x_i)=1$.

Now for $1\leq i,j\leq n$, define the $\Z$--linear map
$\a_{ij}\colon\thinspace \Y_n\rightarrow\A_n$ by its action on
monomials in $\Y_n$:
\[
\a_{ij}(y_{i_1}y_{i_2}\cdots y_{i_{m-1}}y_{i_m}) = a_{ii_1}
a_{i_1i_2} \cdots a_{i_{m-1}i_m} a_{i_m j}
\]
It is then easy to check that $\a_{ij}$ descends to a map on
$\Y_n/\langle y_1^2+2y_1,\ldots,y_n^2+2y_n\rangle$. Finally, if
$\gamma(0)=p_i$ and $\gamma(1)=p_j$, then we set $\psi(\gamma) =
\a_{ij}\circ\y\circ\x(\gamma)$.

\begin{figure}[ht!]\anchor{fig:cords}
\centerline{\small
\psfrag {g}{$\gamma$}
\psfrag {g31}{$\gamma_{31}$}
\psfrag {p1}{$p_1$}
\psfrag {p2}{$p_2$}
\psfrag {p3}{$p_3$}
\includegraphics[width=4in]{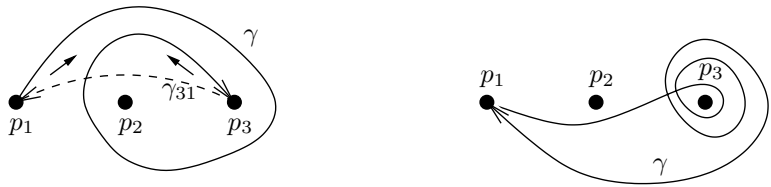}
}
\caption{Cords in $\tilde{\P}_3$}
\label{fig:cords}
\end{figure}

As examples, consider the cords depicted in \figref{fig:cords}.
For the cord $\gamma$ on the left, we can concatenate with
$\gamma_{31}$ and push off of $p_1$ and $p_3$ in the directions
drawn; the resulting loop represents $x_3^{-1}x_2^{-1} \in \F_3$. We
then compute that $\y(\x(\gamma)) = (1+y_3)(1+y_2)$ and
\[
\psi(\gamma) = \a_{13}((1+y_3)(1+y_2)) =
-a_{13}+a_{12}a_{23}+a_{13}a_{32}a_{23}.
\]
This agrees with the definition of $\psi(\gamma)$ from
Section~\ref{sec:braidrep}: since $\gamma =
\sigma_2^{-2}\cdot\gamma_{13}$, we have $\psi(\gamma) =
\hom{\sigma_2}^{-2}(\gamma_{13})$.

For the cord $\gamma$ on the right of \figref{fig:cords}, we
have $\x(\gamma) = x_2x_3^3x_2^{-1}$ and $\y(\x(\gamma))$ $ =
-(1+y_2)(1+y_3)^3(1+y_2) = -(1+y_2)(1+y_3)(1+y_2)$. It follows that
\begin{eqnarray*}
\psi(\gamma) &=& -\a_{11}((1+y_2)(1+y_3)(1+y_2)) \\
&=& 2 - a_{13}a_{31}-a_{12}a_{23}a_{31}
-a_{13}a_{32}a_{21}-a_{12}a_{23}a_{32}a_{21}.
\end{eqnarray*}

\begin{proposition}
\label{prop:tildepsidefined}
$\psi=\alpha\circ\y\circ\x\colon\thinspace
\tilde{\P}_n\rightarrow\A_n$ is well-defined and agrees on $\P_n$
with the definition of $\psi$ from Section~\ref{sec:braidrep}. It
satisfies the skein relation (\ref{eq:psiskein}), even in the case
in which the depicted puncture is an endpoint of the path (so that
there is another component of the path in the depicted neighborhood,
with an endpoint at the puncture).
\end{proposition}

\begin{proof}
To show that $\psi$ is well-defined despite the indeterminacy of
$\gamma$, it suffices to verify that $(\a_{ij}\circ\y)(x_i x) =
(\a_{ij}\circ\y)(x x_j) = (\a_{ij}\circ\y)(x)$ for all $i,j$ and
$x\in\F_n$. This in turn follows from the identity
\begin{eqnarray*}
(\a_{ij}\circ\y)((-1-y_i) y_{i_1}\cdots y_{i_m}) &=&
-a_{ii_1}\cdots a_{i_mj} - a_{ii}a_{ii_1}\cdots a_{i_mj} \\
&=& a_{ii_1}\cdots a_{i_mj} \\
&=&(\a_{ij}\circ\y)(y_{i_1}\cdots y_{i_m})
\end{eqnarray*}
for any $i_1,\ldots,i_m$, and a similar calculation for
$(\a_{ij}\circ\y)(y_{i_1}\cdots y_{i_m}(-1-y_j))$.

We next note that we can set $\x(\gamma_{ij})=1$ by pushing the
relevant loop into the upper half plane; hence $\psi(\gamma_{ij}) =
\a_{ij}(1) = a_{ij}$, which agrees with the normalization from
Proposition~\ref{prop:psi}. Since normalization and the skein
relation (\ref{eq:psiskein}) define $\psi$ on $\P_n$ by
Proposition~\ref{prop:psiskein}, we will be done if we can prove
that the skein relation is satisfied for $\psi =
\a_{ij}\circ\y\circ\x$.

In the skein relation, let $p_k$ be the depicted puncture, and
suppose that the paths on either side of the puncture begin at $p_i$
and end at $p_j$. Then there exist $x,x'\in\F_n$, with $x$ going
from $p_i$ to $p_k$ and $x'$ going from $p_k$ to $p_j$, such that
the two paths avoiding $p_k$ are mapped by $\x$ to $xx'$ and
$xx_kx'$, while the two paths through $p_k$ are mapped to $x$ and
$x'$. The skein relation then becomes
\[
\a_{ij}(\y(xx')) + \a_{ij}(\y(xx_kx')) +
\a_{ik}(\y(x))\a_{kj}(\y(x')) = 0,
\]
which holds by the definitions of $\y$ and $\a_{ij}$:
\[
\a_{ij}(\y(xx')) + \a_{ij}(\y(xx_kx')) = -\a_{ij}(\y(x)x_k\y(x')) =
-\a_{ik}(\y(x))\a_{kj}(\y(x')),
\]
as desired.
\end{proof}

\begin{proposition}
\label{prop:tildepsisuffice}
For $1\leq i\leq n$, let $\gamma_{ii} \in \tilde{\P}_n$ denote the
trivial loop beginning and ending at $p_i$. Then the skein relation
(\ref{eq:psiskein}), and the normalizations
$\psi(\gamma_{ij})=a_{ij}$ for $i\neq j$ and $\psi(\gamma_{ii})=-2$
for all $i$, completely determine the map $\psi$ on $\tilde{\P}_n$.
Furthermore, for $\gamma\in\tilde{\P}_n$ and $B\in B_n$, we have
$\psi(B\cdot\gamma) = \hom{B}(\psi(\gamma))$.
\end{proposition}

\begin{proof}
The normalizations define $\psi(\gamma)$ when $\x(\gamma)=1$, and
the skein relation then allows us to define $\psi(\gamma)$
inductively on the length of the word $\x(\gamma)$, as in the proof
of Proposition~\ref{prop:tildepsidefined}. Note that the given
normalizations are correct because $\x(\gamma_{ii})=1$ and hence
$\psi(\gamma_{ii}) = (\a_{ii}\circ\y)(1) = a_{ii} = -2$.

The proof that $\psi(B\cdot\gamma) = \hom{B}(\psi(\gamma))$
similarly uses induction: it is true when $\x(\gamma)=1$ by
Proposition~\ref{prop:psi} (in particular, it is trivially true if
$\gamma=\gamma_{ii}$), and it is true for general $\gamma$ by
induction, using the skein relation.
\end{proof}

\subsection{Proofs of Theorems~\ref{thm:mainknot} and \ref{thm:mainbraid}}
\label{ssec:thmproofs}

We are now in a position to prove the main results of this paper,
beginning with the identification of braid contact homology with a
cord ring.

\begin{proof}[Proof of Theorem~\ref{thm:mainbraid}]
Let $B\in B_n$, and recall that we embed $B$ in the solid torus
$M=D\times S^1$ in the natural way. If we view $B$ as an element of
the mapping class group of $(D,P)$, then we can write $M$ as
$D\times [0,1]/\sim$, where $D\times\{0\}$ and $D\times\{1\}$ are
identified via the map $B$; the braid then becomes $P\times
[0,1]/\sim$.

Any cord of $B$ in $M$ (in the sense of Definition~\ref{def:cord})
can be lifted to a path in the universal cover $D\times\R$ of $M$,
whence it can be projected to an element of $\tilde{\P}_n$, ie, a
cord in $(D,P)$ (in the sense of Definition~\ref{def:cordplane}).
There is a $\Z$ action on the set of possible lifts, corresponding
in the projection to the action of the map given by $B$. If we
denote by $\tilde{\P}_n/B$ the set of cords in $(D,P)$ modulo the
action of $B$, then it follows that any cord of $B$ in $M$ yields a
well-defined element of $\tilde{\P}_n/B$.

Now for $\gamma\in\tilde{\P}_n$, we have $\psi(B\cdot\gamma) =
\hom{B}(\psi(\gamma))$ by Proposition~\ref{prop:tildepsisuffice};
hence $\psi\colon\thinspace \tilde{\P}_n\rightarrow \A_n$ descends
to a map $\tilde{\P}_n/B\rightarrow\A_n/\im(1-\hom{B}) = HC_0(B)$.
When we compose this with the map from cords of $B$ to
$\tilde{\P}_n/B$, we obtain a map $\A_B \rightarrow HC_0(B)$. This
further descends to a map $\A_B/\I_B \rightarrow HC_0(B)$, since the
skein relations defining $\I_B$ translate to the skein relations
(\ref{eq:psiskein}) in $\tilde{\P}_n$, which are sent to $0$ by
$\psi$ by Proposition~\ref{prop:tildepsidefined}.

It remains to show that the map $\A_B/\I_B \rightarrow HC_0(B)$ is
an isomorphism. It is clearly surjective since any generator
$a_{ij}$ of $\A_n$ is the image of $\gamma_{ij}$, viewed as a cord
of $B$ via the inclusion $(D,P) = (D\times\{0\},P\times\{0\})
\hookrightarrow (M,B)$. To establish injectivity, we first note that
homotopic cords of $B$ in $M$ are mapped to the same element of
$\tilde{\P}_n/B$, and hence the map $\A_B \rightarrow HC_0(B)$ is
injective. Furthermore, if two elements of $\A_B$ are related by a
series of skein relations, then since $\psi$ preserves skein
relations, they are mapped to the same element of $HC_0(B)$; hence
the quotient map on $\A_B/\I_B$ is injective, as desired.
\end{proof}

\begin{proof}[Proof of Theorem~\ref{thm:mainknot}]
Let the knot $K$ be the closure of a braid $B\in B_n$; we picture
$B$ inside a solid torus $M$ as in the above proof, and then embed
(the interior of) $M$ in $\R^3$ as the complement of some line
$\ell$. The braid $B$ in $\R^3\setminus \ell$ thus becomes the knot
$K$ in $\R^3$. It follows that $\A_K/\I_K$ is simply a quotient of
$\A_B/\I_B$, where we mod out by homotopies of cords which pass
through $\ell$.

A homotopy passing through $\ell$ simply replaces a cord of the form
$\gamma_1\gamma_2$ with a cord of the form
$\gamma_1\gamma_*\gamma_2$, where $\gamma_1$ begins on $K$ and ends
at some point $p\in\R^3\setminus\ell$ near $\ell$, $\gamma_2$ begins
at $p$ and ends on $K$, and $\gamma_*$ is a loop with base point $p$
which winds around $\ell$ once. Furthermore, we may choose a point
$*\in D$, near the boundary, such that $p$ corresponds to $(*,0)\in
D\times S^1$, and $\gamma_*$ corresponds to $\{*\}\times S^1$.

The homeomorphism of $(D,P)$ given by $B$ induces a foliation on the
solid torus: if we identify the solid torus with $D\times
[0,1]/\sim$ as in the proof of Theorem~\ref{thm:mainbraid}, then the
leaves of the foliation are given locally by $\{q\}\times [0,1]$ for
$q\in D$. In addition, since $*$ is near the boundary, it is
unchanged by $B$, and so $\gamma_*$ is a leaf of the foliation. We
can use the foliation to project $\gamma_1,\gamma_2$ to cords in
$(D,P\cup\{*\})$, where $D$ is viewed as $D\times\{0\} \subset
D\times S^1$. (This is precisely the projection used in the proof of
Theorem~\ref{thm:mainbraid}.) If we write $\tilde{P}_*^1$ (resp.\
$\tilde{\P}_*^2$) as the set of cords in $(D,P\cup\{*\})$ ending
(resp.\ beginning) at $*$, then $\gamma_1,\gamma_2$ project to cords
$\gamma_1'\in\tilde{P}_*^1,\gamma_2'\in\tilde{P}_*^2$.

Under this projection, the homotopy passing through $\ell$ replaces
the cord $\overline{\gamma_1'\gamma_2'}$ in $(D,P)$ with the cord
$\overline{(\gamma_1')(B\cdot\gamma_2')}$, where
$\overline{\gamma_1'\gamma_2'}$ denotes the cord given by
concatenating the paths $\gamma_1'$ and $\gamma_2'$, and so forth.
To compute $\A_K/\I_K$ from $\A_B/\I_B$, we need to mod out by the
relation which identifies these two cords, for any choice of
$\gamma_1'\in\tilde{\P}_*^1$ and $\gamma_2'\in\tilde{\P}_*^2$. By
using the skein relations in $\A_B/\I_B$, it suffices to consider
the case where $\gamma_1'=\gamma_{i*}$ and $\gamma_2'=\gamma_{*j}$
for some $i,j$, with notation as in Section~\ref{sec:braidrep}. In
this case, we have $\overline{\gamma_1'\gamma_2'}$ homotopic to
$\gamma_{ij}$, while
\[
\psi(\overline{(\gamma_1')(B\cdot \gamma_2')}) = \sum_k a_{ik}
(\Phir{B})_{kj}
\]
by the definition of $\Phir{}$.

It follows that $\A_K/\I_K = \A_n/I$, where $I$ is generated by the
image of $1-\hom{B}$ and by the entries of the matrix
$A-A\cdot\Phir{B}$. Now by Proposition~\ref{prop:matrix}, we have
the matrix identity
\[
((1-\hom{B})(a_{ij})) = A-\Phil{B}\cdot A\cdot\Phir{B} =
(A-\Phil{B}\cdot A) + \Phil{B}\cdot(A-A\cdot\Phir{B}),
\]
and so $I$ is also generated by the entries of the matrices
$A-\Phil{B}\cdot A$ and $A-A\cdot\Phir{B}$.
\end{proof}

\section{Methods to calculate the cord ring}
\label{sec:plat}

So far, we have given only one way to compute the cord ring of a
knot: express the knot as the closure of a braid, and then compute
$HC_0(K)$ using Definition~\ref{def:knotcordring}. In many
circumstances, it is easier to use alternative methods. In this
section, we discuss two such methods. The first technique relies on
a plat presentation of the knot; we describe how to calculate the
cord ring from a plat in Section~\ref{ssec:plat}. We apply this in
Section~\ref{ssec:2bridge} to the case of general two-bridge knots,
for which the cord ring can be explicitly computed in terms of the
determinant. In Section~\ref{ssec:diagram}, we present another
method for calculating the cord ring, this time in terms of any knot
diagram.

\subsection{The cord ring in terms of plats}
\label{ssec:plat}

In this section, we express the cord ring for a knot $K$ in terms of
a plat presentation of $K$. We assume throughout the section that
$K$ is the plat closure of a braid $B\in B_{2n}$; that is, it is
obtained from $B$ by joining together strands $2i-1$ and $2i$ on
each end of the braid, for $1 \leq i\leq n$.

Let $\I_B^{\plat}\subset\A_{2n}$ be the ideal generated by
$a_{ij}-a_{i'j'}$ and $\hom{B}(a_{ij})-\hom{B}(a_{i'j'})$, where
$i,j,i',j'$ range over all values between $1$ and $2n$ inclusive
such that $\lceil i/2 \rceil = \lceil i'/2 \rceil$ and $\lceil j/2
\rceil = \lceil j'/2 \rceil$.

\begin{theorem}
\label{thm:plat}
If $K$ is the plat closure of $B\in B_{2n}$, then the cord ring of
$K$ is isomorphic to $\A_{2n}/\I_B^{\plat}$.
\end{theorem}

Note that $\A_{2n}/\I_B^{\plat}$ can be expressed as a quotient of
$\A_n$, as follows. Define an algebra map $\eta\colon\thinspace
\A_{2n}\rightarrow\A_n$ by $\eta(a_{ij}) = a_{\lceil i/2 \rceil,
\lceil j/2 \rceil}$. Then $\eta$ induces an isomorphism
$\A_{2n}/\I_B^{\plat} \cong \A_n/\eta(\hom{B}(\ker\eta))$, where
$\eta(\hom{B}(\ker\eta))$ is the ideal in $\A_n$ given by the image
of $\ker\eta\subset\A_{2n}$ under the map $\eta\circ\hom{B}$.

Calculating the cord ring using Theorem~\ref{thm:plat} is reasonably
simple for small knots. For example, the trefoil is the plat closure
of $\sigma_2^3\in B_4$. Here are generators of the kernel of
$\eta\,:\A_4\rightarrow\A_2$, along with their images under
$\eta\circ\hom{\sigma_2^3}$:
\begin{gather*}
\begin{align*}{2}
2+a_{12} &\mapsto 2-3a_{12}+a_{12}a_{21}a_{12} &
2+a_{34} &\mapsto 2-3a_{12}+a_{12}a_{21}a_{12} \\
2+a_{21} &\mapsto 2-3a_{21}+a_{21}a_{12}a_{21} &
2+a_{43} &\mapsto 2-3a_{21}+a_{21}a_{12}a_{21} \\
a_{14}-a_{13} &\mapsto -2+a_{12}+a_{12}a_{21} &
a_{41}-a_{31} &\mapsto -2+a_{21}+a_{12}a_{21} \\
a_{14}-a_{23} &\mapsto a_{12}-a_{21} & a_{41}-a_{32} &\mapsto
a_{21}-a_{12}
\end{align*} \\
\begin{align*}
a_{14}-a_{24} &\mapsto
2+a_{12}-4a_{21}a_{12}+a_{21}a_{12}a_{21}a_{12} \\
a_{41}-a_{42} &\mapsto
2+a_{21}-4a_{21}a_{12}+a_{21}a_{12}a_{21}a_{12}
\end{align*}
\end{gather*}
Since $a_{12}-a_{21}\in\eta(\hom{\sigma_2^3}(\ker\eta))$,
we set $x:=-a_{12}=-a_{21}$ in
$\A_2/\eta(\hom{\sigma_2^3}(\ker\eta))$. The above images then give
the relations $2+3x-x^3$, $-2-x+x^2$, $2-x-4x^2+x^4$, with gcd
$-2-x+x^2$, and so $HC_0(3_1) \cong \Z[x]/(x^2-x-2)$.

\begin{proof}[Proof of Theorem~\ref{thm:plat}]
Embed $B\in B_{2n}$ in $D\times [0,1]$, so that the endpoints of $B$
are given by $(p_i,0)$ and $(p_i,1)$ for $1\leq i\leq 2n$ and some
points $p_1,\ldots,p_{2n}\in D$. (See \figref{fig:platex} for an
example.) As in the proof of Theorem~\ref{thm:mainbraid}, any cord
of $B$ in $D\times [0,1]$ can be isotoped to a cord of $(D,P)$,
where $P=\{p_1,\ldots,p_{2n}\}$ and $(D,P)$ is viewed as
$(D\times\{0\},P\times\{0\}) \subset (D\times [0,1],B)$. Hence there
is a map from cords of $B$ to $\A_{2n}$ induced by the map $\psi$
from Section~\ref{ssec:cordmap}, and this map respects the skein
relations (\ref{eq:skein1}), (\ref{eq:skein2}).

\begin{figure}[ht!]\anchor{fig:platex}
\centerline{
\includegraphics[width=3in]{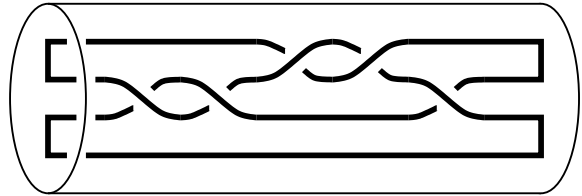}
}
\caption{
Plat representation of the knot $5_2$ in $D\times [0,1]$
}
\label{fig:platex}
\end{figure}

We may assume that $p_1,\ldots,p_{2n}$ lie in order on a line in
$D$; then $K$ is the union of $B\subset D\times [0,1]$ and the line
segments $L_j\times\{0\}$ and $L_j\times\{1\}$, where $1\leq j\leq
n$ and $L_j$ connects $p_{2j-1}$ and $p_{2j}$. Any cord of $K$ in
$\R^3$ can be isotoped to a cord lying in $D \times (0,1)$, by
``pushing'' any section lying in $\R^2\times (-\infty,0]$ or
$\R^2\times [1,\infty)$ into $\R^2\times (0,1)$, and then
contracting $\R^2$ to $D$. To each cord $\gamma$ of $K$, we can thus
associate a (not necessarily unique) element of $\A_{2n}$, which we
denote by $\psi(\gamma)$.

Because of the line segments $L_j\times\{0\}$, isotopic cords of $K$
may be mapped to different elements of $\A_{2n}$. More precisely,
any cord with an endpoint at $(p_{2j-1},0)$ is isotopic via
$L_j\times\{0\}$ to a corresponding cord with endpoint at
$(p_{2j},0)$, and vice versa. To mod out by these isotopies, we mod
out $\A_{2n}$ by $a_{ij}-a_{i'j'}$ for all $i,j,i',j'$ with $\lceil
i/2 \rceil = \lceil i'/2 \rceil$ and $\lceil j/2 \rceil = \lceil
j'/2 \rceil$. Similarly, isotopies using the line segments
$L_j\times\{1\}$ require that we further mod out $\A_{2n}$ by
$\hom{B}(a_{ij})-\hom{B}(a_{i'j'})$ for the same $i,j,i',j'$; note
that $\hom{B}$ appears because all cords must be translated from
$D\times\{1\}$ to $D\times\{0\}$.

\begin{figure}[ht!]\anchor{fig:slip}
\centerline{
\includegraphics[width=1.6in]{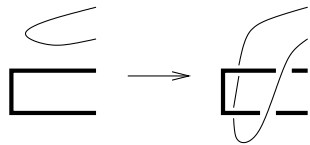}
}
\caption{
Slipping a segment of a cord around $L_j\times\{0\}$
}
\label{fig:slip}
\end{figure}

We now have a map, which we also write as $\psi$, from cords of $K$
to $\A_{2n}/\I_B^{\plat}$, which satisfies the skein relations
(\ref{eq:skein1}), (\ref{eq:skein2}). (Any skein relation involving
one of the segments $L_j\times\{0\}$ or $L_j\times\{1\}$ can be
isotoped to one which involves a section of $B$ instead.) To ensure
that the map is well-defined, we still need to check that the
particular isotopy from a cord of $K$ to a cord in $D\times (0,1)$
is irrelevant. That is, the isotopy shown in \figref{fig:slip}
should not affect the value of $\psi$. (There is a similar isotopy
around $L_j\times\{1\}$ instead of $L_j\times\{0\}$, which can be
dealt with similarly.)

In the projection to $(D,P)$, the isotopy pictured in
\figref{fig:slip} corresponds to moving a segment of a cord on
one side of $L_j$ across to the other side, ie, passing this
segment through the points $p_{2j-1}$ and $p_{2j}$. Now we have the
following chain of equalities in $\A_{2n}/\I_B^{\plat}$:
\begin{eqnarray*}
\psi(\raisebox{-0.26in}{\includegraphics[height=0.6in]{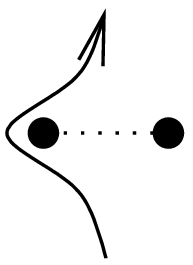}})
&=&
-\psi(\raisebox{-0.26in}{\includegraphics[height=0.6in]{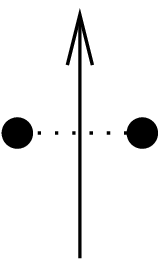}})
-\psi(\raisebox{-0.26in}{\includegraphics[height=0.6in]{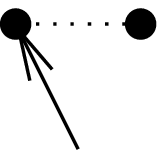}})
\psi(\raisebox{-0.26in}{\includegraphics[height=0.6in]{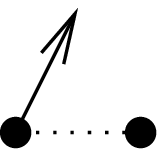}}) \\
&=&
\psi(\raisebox{-0.26in}{\includegraphics[height=0.6in]{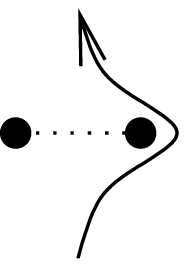}})
-\psi(\raisebox{-0.26in}{\includegraphics[height=0.6in]{figures/loop3.eps}})
\psi(\raisebox{-0.26in}{\includegraphics[height=0.6in]{figures/loop4.eps}})
+\psi(\raisebox{-0.26in}{\includegraphics[height=0.6in]{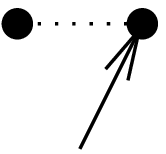}})
\psi(\raisebox{-0.26in}{\includegraphics[height=0.6in]{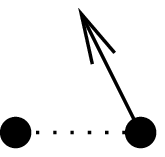}}) \\
&=&
\psi(\raisebox{-0.26in}{\includegraphics[height=0.6in]{figures/loop5.eps}})
\end{eqnarray*}
where the dotted line represents $L_j$, and the last equality holds
by the definition of $\I_B^{\plat}$. Hence the value of $\psi$ is
unchanged under the isotopy of \figref{fig:slip}.

To summarize, we have a map from cords of $K$, modulo homotopy and
skein relations, to $\A_{2n}/\I_B^{\plat}$. By construction, this
induces an isomorphism between the cord ring of $K$ and
$\A_{2n}/\I_B^{\plat}$, as desired.
\end{proof}

The argument of Theorem~\ref{thm:plat} also gives a plat description
of $HC_0^{\ab}(K)$.

\begin{corollary}
\label{cor:platab}
$HC_0^{\ab}(K)$ can be obtained from $HC_0(K) \cong
\A_{2n}/\I_B^{\plat}$ by further quotienting by $a_{ij}-a_{ji}$ for
all $i,j$ and abelianizing. The result can be viewed as a quotient
of the polynomial ring $\Z[\{a_{ij}|1\leq i<j\leq n\}]$.
\end{corollary}

\subsection{Two-bridge knots}
\label{ssec:2bridge}

For two-bridge knots, Theorem~\ref{thm:plat} implies that $HC_0$ has
a particularly simple form. In particular, $HC_0$ is a quotient of
$\A_2 = \Z\langle a_{12},a_{21}\rangle$; in this quotient, it turns
out that $a_{12}=a_{21}$, so that $HC_0$ is a quotient of a
polynomial ring $\Z[x]$. The main result of this section shows that
for two-bridge knots, $HC_0$ is actually determined by the knot's
determinant.

We recall some notation from \cite{I}. If $K$ is a knot,
$\Delta_K(t)$ denotes the Alexander polynomial of $K$ as usual, and
$|\Delta_K(-1)|$ is the determinant of $K$. Define the sequence of
polynomials $\{p_m\in\Z[x]\}$ by $p_0(x)=2-x$, $p_1(x)=x-2$,
$p_{m+1}(x)=xp_m(x)-p_{m-1}(x)$.

\begin{theorem}
\label{thm:2bridge}
If $K$ is a $2$--bridge knot, then
\[
HC_0(K) \cong HC_0^{\ab}(K) \cong
\Z[x]/(p_{(|\Delta_K(-1)|+1)/2}(x)).
\]
\end{theorem}

\noindent This generalizes \cite[Proposition 7.3]{I}. Also compare this
result to \cite[Theorem 6.13]{I}, which states that for any knot $K$,
there is a surjection from $HC_0(K)$ to $\Z[x]/(p_{(n(K)+1)/2}(x))$,
where $n(K)$ is the largest invariant factor of the first homology
of the double branched cover of $K$.

Before we can prove Theorem~\ref{thm:2bridge}, we need to recall
some results (and more notation) from \cite{I}, and establish a few
more lemmas. Define the sequence $\{q_m\in\Z[x]\}$ by $q_0(x)=-2$,
$q_1(x)=-x$, $q_{m+1}(x) = xq_m(x)-q_{m-1}(x)$; this recursion
actually defines $q_m$ for all $m\in\Z$, and $q_{-m}=q_m$. The Burau
representation of $B_n$ with $t=-1$ is given as follows:
$\Bur{\sigma_k}$ is the linear map on $\Z^n$ whose matrix is the
identity, except for the $2\times 2$ submatrix formed by the $k,k+1$
rows and columns, which is $\left( \begin{smallmatrix} 2 & 1 \\ 1 &
0 \end{smallmatrix} \right)$. This extends to a representation which
sends $B\in B_n$ to $\Bur{B}$. If $B$ is a braid, let $\hat{B}$
denote the braid obtained by reversing the word which gives $B$.

\begin{lemma}\label{lem:Buraumatrix}
{\rm\cite{I}}\qua
For $B\in B_n$ and $v\in\Z^n$, if we set $a_{ij} = q_{v_i-w_j}$ for
all $i,j$, then $\hom{B}(a_{ij}) =
q_{(\Bur{\hat{B}}v)_i-(\Bur{\hat{B}}v)_j}$ for all $i,j$.
\end{lemma}

Let $K$ be a $2$--bridge knot; then $K$ is the plat closure of some
braid in $B_4$ of the form
$B=\sigma_2^{-a_1}\sigma_1^{b_1}\sigma_2^{-a_2}\cdots
\sigma_1^{b_k}\sigma_2^{-a_{k+1}}$. As usual, we can then associate
to $K$ the continued fraction
\[
\frac{m}{n} = a_1 + \frac{1}{b_1} \!\!
\begin{array}{l}
\\ + \end{array}
\!\! \frac{1}{a_2} \!\!
\begin{array}{l}
\\ + \, \cdots \, + \end{array}
\!\! \frac{1}{b_k}
\begin{array}{l}
\\ + \end{array}
\!\! \frac{1}{a_{k+1}}\, ,
\]
where $\gcd(m,n)=1$ and $n>0$.

\begin{lemma}
\label{lem:Burauvector}
For $B,m,n$ as above, we have $\Bur{B} \left( \begin{smallmatrix} 0
\\ 0 \\ 1 \\ 1 \end{smallmatrix} \right) = \left(
\begin{smallmatrix} -n+1 \\ m-n+1 \\ m+1 \\ 1 \end{smallmatrix}
\right).$
\end{lemma}

\begin{proof}
We can compute that
\[
\textstyle{ \Bur{\sigma_1} \left( \begin{smallmatrix} -n+1 \\ m-n+1
\\ m+1 \\ 1 \end{smallmatrix} \right) = \left( \begin{smallmatrix}
-m-n+1 \\ -n+1 \\ m+1 \\ 1 \end{smallmatrix} \right)
}
\hspace{0.15in} \textrm{and} \hspace{0.15in} \textstyle{
\Bur{\sigma_2^{-1}} \left( \begin{smallmatrix} -n+1 \\ m-n+1 \\ m+1
\\ 1 \end{smallmatrix} \right) = \left( \begin{smallmatrix} -n+1 \\
m+1 \\ m+n+1 \\ 1 \end{smallmatrix} \right).
}
\]
The lemma follows easily by induction.
\end{proof}

We present one final lemma about the polynomials $p_k$ and $q_k$.
For any $k,m$, define $r_{k,m} = q_k-q_{k-m}$.

\begin{lemma}
\label{lem:gcd}
If $m>0$ is odd and $\gcd(m,n)=1$, then $\gcd(r_{m,m},r_{n,m}) =
p_{(m+1)/2}$.
\end{lemma}

\begin{proof}
We first note that $r_{m-k,m}=-r_{k,m}$ and
$r_{2k-l,m}-r_{l,m}-q_{l-k} r_{k,m} = 0$ for all $k,l,m$; the first
identity is obvious, while the second is an easy induction (or use
Lemma 6.14 from \cite{I}). It follows that $\gcd(r_{k,m},r_{l,m})$
is unchanged if we replace $(k,l)$ by any of $(l,k)$, $(k,m-l)$,
$(l,2k-l)$.

Now consider the operation which replaces any ordered pair $(k,l)$
for $k,l>m/2$ with: $(l,2l-k)$ if $k\geq l$ and $2l-k>m/2$;
$(l,m-2l+k)$ if $k\geq l$ and $2l-k<m/2$; $(k,2k-l)$ if $k<l$ and
$2k-l>m/2$; $(k,m-2k+l)$ if $k<l$ and $2k-l<m/2$. This operation
preserves $\gcd(r_{k,m},r_{l,m})$ and $\gcd(2k-m,2l-m)$, as well as
the condition $k,l>m/2$; it also strictly decreases $\max(k,l)$
unless $k=l$.

We can now use a descent argument, beginning with the ordered pair
$(m,n)$ and performing the operation repeatedly until we obtain a
pair of the form $(k,k)$. We then have $2k-m=\gcd(2m-n,m)=1$ and
$\gcd(r_{m,m},r_{n,m}) = r_{k,m}$. The lemma now follows from the
fact, established by induction, that $r_{(m+1)/2,m} = p_{(m+1)/2}$.
\end{proof}

\begin{proof}[Proof of Theorem~\ref{thm:2bridge}]
As usual, we assume that $K$ is the plat closure of $B\in B_4$; it
is then also the case that $K$ is the plat closure of $\hat{B}$. Let
$m/n$ be the continued fraction associated to $K$, and note that
$|\Delta_K(-1)|=m$.

By Theorem~\ref{thm:plat}, $HC_0(K)$ is a quotient of $\A_4$, and in
this quotient, we can set $a_{12}=a_{21}=a_{34}=a_{43}=-2$,
$a_{13}=a_{14}=a_{23}=a_{24}$, and $a_{31}=a_{41}=a_{32}=a_{42}$.

We first compute $HC_0^{\ab}(K)$. Here we can further set
$a_{13}=a_{31}=:-x$. Then we have $a_{ij} = q_{v_i-v_j}$, where $v$
is the vector $(0,0,1,1)$. By Lemma~\ref{lem:Buraumatrix}, we have
$\hom{\hat{B}}(a_{ij}) = q_{(\Bur{B}v)_i-(\Bur{B}v)_j}$; by
Lemma~\ref{lem:Burauvector}, we conclude the matrix identity:
\[
(\hom{\hat{B}}(a_{ij})) = \left( \begin{array}{cc|cc}
q_0 & q_{-m} & q_{-m-n} & q_{-n} \\
q_m & q_0 & q_{-n} & q_{m-n} \\ \hline
q_{m+n} & q_n & q_0 & q_m \\
q_n & q_{-m+n} & q_{-m} & q_0 \end{array} \right) 
\]
Here we have divided the matrix into $2\times 2$ blocks for clarity.
By Theorem~\ref{thm:plat}, the relations defining $HC_0(K)$ then
correspond to equating the entries within each block. In other
words, since $q_{-k}=q_k$ for all $k$, we have
\[
HC_0(K) \cong \Z[x]/(q_m-q_0,q_{m+n}-q_n,q_n-q_{n-m}) =
\Z[x]/(r_{m,m},r_{m+n,m},r_{n,m}).
\]
We may assume without loss of generality that $m>0$ (otherwise
replace $m$ by $-m$); furthermore, $m$ is odd since $K$ is a knot
rather than a two-component link. By Lemma~\ref{lem:gcd}, we can
then conclude that $HC_0(K) \cong \Z[x]/(p_{(m+1)/2}) =
\Z[x]/(p_{(|\Delta_K(-1)|+1)/2})$.

The computation of $HC_0(K)$ rather than $HC_0^{\ab}(K)$ is very
similar but becomes notationally more complicated. Set $a_{13}=a_1$
and $a_{31}=a_2$, so that $HC_0(K)$ is a quotient of $\Z\langle
a_1,a_2 \rangle$. As in Section 7.3 of \cite{I}, we define two
sequences $\{q_k^{(1)},q_k^{(2)}\}$ by $q_0^{(1)}=q_0^{(2)}=-2$,
$q_1^{(1)}=a_1$, $q_1^{(2)}=a_2$, and $q_{m+1}^{(1)} = -a_1
q_m^{(2)}-q_{m-1}^{(1)}$, $q_{m+1}^{(2)} = -a_2
q_m^{(1)}-q_{m-1}^{(2)}$. Note that $q_m^{(1)}|_{a_1=a_2=-x} =
q_m^{(2)}|_{a_1=a_2=-x} = q_m$, and that each nonconstant monomial
in $q_m^{(1)}$ (resp.\ $q_m^{(2)}$) begins with $a_1$ (resp.\
$a_2$).

In terms of $a_1,a_2$, the monomials appearing in
$\hom{\hat{B}}(a_{ij})$ look like $a_1a_2a_1\cdots$ or
$a_2a_1a_2\cdots$. Since $\hom{\hat{B}}(a_{ij})$ projects to the
appropriate polynomial $q_k$ if we set $a_1=a_2=-x$, it readily
follows that
\[
(\hom{\hat{B}}(a_{ij})) = \left( \begin{array}{cc|cc}
q_0^{(r)} & q_{-m}^{(r)} & q_{-m-n}^{(r)} & q_{-n}^{(r)} \\
q_m^{(s)} & q_0^{(s)} & q_{-n}^{(s)} & q_{m-n}^{(s)} \\ \hline
q_{m+n}^{(1)} & q_n^{(1)} & q_0^{(1)} & q_m^{(1)} \\
q_n^{(2)} & q_{-m+n}^{(2)} & q_{-m}^{(2)} & q_0^{(2)} \end{array}
\right)
\]
where $(r,s)=(1,2)$ or $(2,1)$. (The superscripts follow from an
inspection of the permutation on the four strands induced by
$\hat{B}$.) To obtain $HC_0(K)$ from $\Z\langle a_1,a_2 \rangle$, we
quotient by setting the entries of each $2\times 2$ block equal to
each other.

If we define $r_{k,m}^{(1)}=q_k^{(1)}-q_{k-m}^{(1)}$,
$r_{k,m}^{(2)}=q_k^{(2)}-q_{k-m}^{(2)}$,
$s_{k,m}^{(1)}=q_k^{(1)}-q_{k-m}^{(2)}$,
$s_{k,m}^{(2)}=q_k^{(2)}-q_{k-m}^{(1)}$, then we have
\[
r_{2l-k,m}^{(1)}-r_{k,m}^{(1)}-q_{k-l}^{(p)}r_{l,m}^{(1)}=0
\]
for all $k,l,m$, where $p=1$ or $2$ depending on the parity of
$k-l$, with similar relations for $r^{(2)}$, $s^{(1)}$, and
$s^{(2)}$. Using these identities and the descent argument of
Lemma~\ref{lem:gcd}, we deduce (after a bit of work) that
\[
HC_0(K) \cong \Z\langle a_1,a_2\rangle / \langle
r_{(m+1)/2,m}^{(1)}, r_{(m+1)/2,m}^{(2)}, s_{(m+1)/2,m}^{(1)},
s_{(m+1)/2,m}^{(2)} \rangle.
\]
It can be directly deduced at this point that $a_1=a_2$ in
$HC_0(K)$, whence we can argue as before, but we can circumvent this
somewhat involved calculation by noting that we have now established
that $HC_0(K)$ depends only on $m=|\Delta_K(-1)|$. Since $T(2,2m-1)$
is a $2$--bridge knot with determinant $m$, it follows that $HC_0(K)$
is isomorphic to $HC_0(T(2,2m-1))$, which is $\Z[x]/(p_{(m+1)/2})$
by \cite[Proposition 7.2]{I}.
\end{proof}

We conclude this section by noting that Theorem~\ref{thm:plat} and
Corollary~\ref{cor:platab} are also useful for knots that are not
$2$--bridge. For instance, if $K$ has bridge number $3$, then
$HC_0^{\ab}(K)$ is a quotient of $\Z[a_{12},a_{13},a_{23}]$, and can
be readily computed in many examples, given Gr\"obner basis software
and sufficient computer time. One can calculate, for example, that
$HC_0^{\ab}(P(3,3,2)) \cong \Z[x]/((x-1)p_{11})$ and
$HC_0^{\ab}(P(3,3,-2)) \cong \Z[x]/((x-1)p_5)$, where $P(p,q,r)$ is
the $(p,q,r)$ pretzel knot. For other knots, such as $P(3,3,3)$ and
$P(3,3,-3)$, $HC_0^{\ab}$ is not a quotient of $\Z[x]$. See also
Section~\ref{ssec:sigma2}.

\subsection{The cord ring in terms of a knot diagram}
\label{ssec:diagram}

Here we give a description of the cord ring of a knot given any knot
diagram, not necessarily a plat or a braid closure.

Suppose that we are given a knot diagram for $K$ with $n$ crossings.
There are $n$ components of the knot diagram (ie, segments of $K$
between consecutive undercrossings), which we may label
$1,\ldots,n$. For any $i,j$ in $\{1,\ldots,n\}$, we can define a
cord $\gamma_{ij}$ of $K$ which begins at any point on component
$i$, ends at any point on component $j$, and otherwise lies
completely above the plane of the knot diagram. (In particular, away
from a neighborhood of each endpoint, it lies above any crossings of
the knot.) Such a cord is well-defined up to homotopy.

\begin{figure}[ht!]\anchor{fig:diagram}
\centerline{
\small
\psfrag {glj}{$\gamma_{lj}$}
\psfrag {glk}{$\gamma_{lk}$}
\psfraga <-3pt,0pt> {gli}{$\gamma_{li}$}
\psfraga <-2pt,0pt> {i}{$i$}
\psfrag {j}{$j$}
\psfrag {k}{$k$}
\psfrag {l}{$l$}
\includegraphics[height=1.2in]{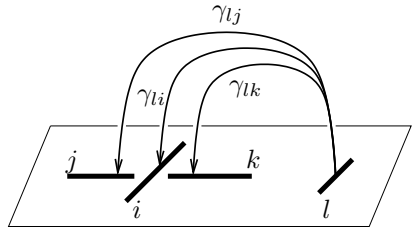}
}
\caption{
The cords $\gamma_{li}$, $\gamma_{lj}$, $\gamma_{lk}$ are related by
a skein relation.
}
\label{fig:diagram}
\end{figure}

Any cord of $K$ can be expressed, via the skein relations, in terms
of these cords $\gamma_{ij}$; imagine pushing the cord upwards while
fixing its endpoints, using the skein relations if necessary, until
the result consists of cords which lie completely above the plane of
the diagram. The crossings in the knot diagram give relations in the
cord ring. More precisely, consider a crossing whose overcrossing
strand is component $i$, and whose undercrossing strands are
components $j$ and $k$. For any $l$, the cords $\gamma_{lj}$ and
$\gamma_{lk}$ are obtained from one another by passing through
component $i$; since the cord joining overcrossing to undercrossing
is $\gamma_{ij}$ (which is homotopic to $\gamma_{ik}$), we have the
skein relation $\gamma_{lj} + \gamma_{lk} + \gamma_{li} \cdot
\gamma_{ij} = 0$. See \figref{fig:diagram}. Similarly, there are
skein relations of the form $\gamma_{jl} + \gamma_{kl} + \gamma_{ji}
\cdot \gamma_{il} = 0$.

Let $\A_n$ denote the usual tensor algebra, and set $a_{ii} = -2$
for all $i$. Define $\I^{\diag}_K \subset \A_n$ to be the ideal
generated by the elements $a_{lj} + a_{lk} + a_{li} a_{ij}$, $a_{jl}
+ a_{kl} + a_{ji} a_{il}$, where $l=1,\ldots,n$ and $(i,j,k)$ ranges
over all $n$ crossings of the knot diagram; as before, $i$ is the
overcrossing strand and $j,k$ are the undercrossing strands. Then
there is a map from the cord ring of $K$ to $\A_n/\I^{\diag}_K$
given by sending $\gamma_{ij}$ to $a_{ij}$ for all $i,j$. It is
straightforward to check that this map is well-defined and an
isomorphism. In particular, all skein relations in the cord ring
follow from the skein relations mentioned above.

\begin{proposition}
The cord ring of $K$ is isomorphic to $\A_n/\I^{\diag}_K$.
\label{prop:diagram}
\end{proposition}

This result may seem impractical, because it expresses the cord ring
of a knot with $n$ crossings as a ring with $n(n-1)$ generators.
However, each relation generating $\I^{\diag}_K$ allows us to
express one generator in terms of three others, and this helps in
general to eliminate the vast majority of these generators.

As an example, consider the usual diagram for a trefoil, and label
the diagram components $1,2,3$ in any order. The crossing where $1$
is the overcrossing strand yields relations
$-2+a_{13}+a_{12}a_{23}$, $a_{21}-a_{23}$, $a_{31}-2+a_{32}a_{23}$,
$-2+a_{31}+a_{32}a_{21}$, $a_{12}-a_{32}$, $a_{13}-2+a_{32}a_{23}$.
The other two crossings yield the same relations, but with indices
cyclically permuted. In $\A_n/\I^{\diag}_K$, we conclude that
$a_{12}=a_{32}=a_{31}=a_{21}=a_{23}=a_{13}$, and the cord ring for
the trefoil is $\Z[x]/(x^2+x-2)$.

Among the three techniques we have discussed to calculate the cord
ring (braid closure, plat, diagram), there are instances when the
diagram technique is computationally easiest. In particular, using
diagrams allows us to study the cord ring for a knot in terms of
tangles contained in the knot.

\section{Some geometric remarks}
\label{sec:remarks}

In this section, we discuss some geometric consequences of the cord
ring construction. Section~\ref{ssec:minchords} relates the cord
ring to binormal chords of a knot; Section~\ref{ssec:sigma2}
establishes a close connection between the abelian cord ring and the
double branched cover of the knot; and Section~\ref{ssec:extensions}
discusses some ways to extend the cord ring to other invariants.

\subsection{Minimal chords}
\label{ssec:minchords}

Here we apply the cord ring to deduce a lower bound on the number of
minimal chords (see below for definition) of a knot in terms of the
double branched cover of the knot. We will also indicate a
conjectural way in which the entire knot DGA of \cite{I} could be
defined in terms of chords. Note that the results in this section
are equally valid for links.

If we impose the usual metric on $\R^3$, we can associate a length
to any sufficiently well-behaved ($L^2$) cord of a knot $K$. Define
a \textit{minimal chord}\footnote{Regarding the spelling: following
a suggestion of D Bar-Natan, we have adopted the spelling ``cord''
in this paper so as to avoid confusion with the chords from the
theory of Vassiliev invariants. In this case, however, a minimal
cord is in fact a ``chord,'' that is, a straight line segment.} to
be a nontrivial cord which locally minimizes length. In other words,
the embedding of $K$ in $\R^3$ gives a distance function
$d\colon\thinspace S^1\times S^1 \rightarrow \R_{\geq 0}$, where
$S^1$ parametrizes $K$ and $d$ is the usual distance between two
points in $\R^3$; then a minimal chord is a local minimum for $d$
not lying on the diagonal of $S^1\times S^1$. Clearly any minimal
chord can be traversed in the opposite direction and remains a
minimal chord; when counting minimal chords, we will identify
minimal chords with their opposites and only count one from each
pair. Minimal chords have previously been studied in the literature,
especially in the context of the ``thickness'' or ``ropelength'' of
a knot.

In the cord ring of $K$, any cord can be expressed in terms of
minimal chords. To see this, imagine a cord as a rubber band, and
pull it taut while keeping its endpoints on $K$. If the result is
not a minimal chord, then it ``snags'' on the knot, giving a union
of broken line segments; using the skein relation, we can express
the result in terms of shorter cords, which we similarly pull taut,
and so forth, until all that remains are minimal chords.

\begin{proposition}
The number of minimal chords for any embedding $K\subset\R^3$ is at
least the minimal possible number of generators of the ring
$HC_0^{\ab}(K)$.
\label{prop:minchords}
\end{proposition}

\noindent As a consequence, for instance, any embedding of
$P(3,3,3)$, $P(3,3,-3)$, or $3_1\# 3_1$ in $\R^3$ has at least two
minimal chords.

A similar result which is slightly weaker, but generally easier to
apply, involves the linearized group $HC_0^{\lin}$ introduced in
\cite{I}. We first note the following expression for $HC_0^{\lin}$
in terms of cords, which is an immediate consequence of
Theorem~\ref{thm:mainknot}.

\begin{proposition}
The group $HC_0^{\lin}(K)$ is the free abelian group generated by
homotopy classes of cords of $K$, modulo the relations:
\eject
\[
\raisebox{-0.17in}{\includegraphics[width=0.4in]{figures/skein12.eps}}
+
\raisebox{-0.17in}{\includegraphics[width=0.4in]{figures/skein11.eps}}
- 2
\raisebox{-0.17in}{\includegraphics[width=0.4in]{figures/skein13.eps}}
- 2
\raisebox{-0.17in}{\includegraphics[width=0.4in]{figures/skein14.eps}}
= 0
\]
\[
\raisebox{-0.17in}{\includegraphics[width=0.4in]{figures/skein21.eps}}
= 0
\]
\end{proposition}

\noindent The argument used to establish
Proposition~\ref{prop:minchords} now yields the following.

\begin{proposition}
The number of minimal chords for an embedding of $K$ is at least the
minimal possible number of generators of the group
$HC_0^{\lin,\ab}(K)$.
\label{prop:minchordslin}
\end{proposition}

We can use Proposition~\ref{prop:minchordslin} to give a lower bound
on the number of minimal chords of a knot which involves only
``classical'' topological information, without reference to the cord
ring.

\begin{corollary}
If $K$ is a knot, let $m(K)$ be the number of invariant factors of
the abelian group $H_1(\Sigma_2(K))$, where $\Sigma_2(K)$ is the
double branched cover of $S^3$ over $K$. Then the number of minimal
chords for an embedding of $K$ is at least $\binom{m(K)+1}{2}$.
\label{cor:minchords}
\end{corollary}

\begin{proof}
By \cite[Proposition 7.11]{I} there is a surjection of groups from
$HC_0^{\lin,\ab}(K)$ to $\Sym^2(H_1(\Sigma_2(K)))$. It is easy to
see that the minimal number of generators of
$\Sym^2(H_1(\Sigma_2(K)))$ is $\binom{m(K)+1}{2}$.
\end{proof}

As a result of Corollary~\ref{cor:minchords}, we can demonstrate
that there are knot types for which the number of minimal chords
must be arbitrarily large.

\begin{corollary}
Let $K$ be a knot, and $\overline{K}$ its mirror. The number of
minimal chords of an embedding of the knot $\#^{m_1} K \#^{m_2}
\overline{K}$ is at least $\binom{m_1+m_2+1}{2}$.
\label{cor:minchordconnectsum}
\end{corollary}

\begin{proof}
We have $H_1(\Sigma_2(\#^{m_1} K \#^{m_2} \overline{K})) \cong
\oplus^{m_1+m_2} H_1(\Sigma_2(K))$.
\end{proof}

\noindent It is not hard to show that any knot with bridge number
$k$ has an embedding with exactly $\binom{k}{2}$ minimal chords.
Hence Corollary~\ref{cor:minchordconnectsum} gives a sharp bound
whenever $K$ is $2$--bridge. To the author's knowledge, it is an open
problem to find sharp lower bounds for the number of minimal chords
for a general knot.

The fact that $HC_0(K)$ can be expressed in terms of minimal chords
suggests that there might be a similar expression for the entire
knot DGA (see \cite{I} for definition), of which $HC_0$ is the
degree $0$ homology. Here we sketch a conjectural formulation for
the knot DGA in terms of chords.

Let a \textit{segment chord} of $K$ be a cord consisting of a
directed line segment; note that the space of segment chords is
parametrized by $S^1\times S^1$, minus a $1$--dimensional subset $C$
corresponding to segments which intersect $K$ in an interior point.
Generically, there are finitely many \textit{binormal chords} of
$K$, which are normal to $K$ at both endpoints; these are critical
points of the distance function $d$ on $S^1\times S^1$, and include
minimal chords. The critical points of $d$ then consist of binormal
chords, along with the diagonal in $S^1\times S^1$.

Let $\A$ denote the tensor algebra generated by binormal chords of
$K$, with grading given by setting the degree of a binormal chord to
be the index of the corresponding critical point of $d$. We can
define a differential on $\A$ using gradient flow trees, as we now
explain.

\begin{figure}[ht!]\anchor{fig:bifurcate}
\centerline{
\includegraphics[width=2.8in]{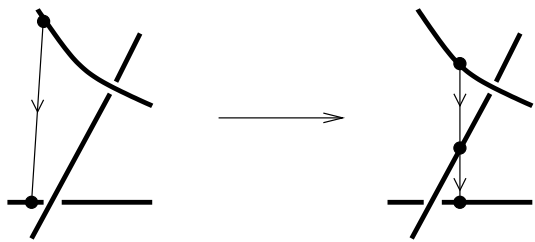}
}
\caption{
Bifurcation in gradient flow. The segment chord on the left can
split into the two chords on the right, each of which subsequently
follows negative gradient flow.
}
\label{fig:bifurcate}
\end{figure}

In the present context, a \textit{gradient flow tree} consists of
negative gradient flow for $d$ on $S^1\times S^1$, except that the
flow is allowed to bifurcate at a point of $C$, by jumping from this
point $(t_1,t_2)$ to the two points $(t_1,t_3),(t_3,t_2)$
corresponding to the segment chords into which $K$ divides the chord
$(t_1,t_2)$. (See \figref{fig:bifurcate}.) Now consider binormal
chords $a_i,a_{j_1},\ldots,a_{j_k}$ (not necessarily distinct), and
look at the moduli space $\M(a_i;a_{j_1},\ldots,a_{j_k})$ of
gradient flow trees beginning at $a_i$ and ending at
$a_{j_1},\ldots,a_{j_k}$, possibly along with some ending points on
the diagonal of $S^1\times S^1$. To each such tree, we associate the
monomial $(-2)^p a_{j_1}\cdots a_{j_k}$, where $p$ is the number of
endpoints on the diagonal, and the order of the $a_{j_l}$'s is
determined in a natural way by the bifurcations. The expected
dimension of $\M(a_i;a_{j_1},\ldots,a_{j_k})$ turns out to be $\deg
a_i - \sum \deg a_{j_l}-1$; we then define the differential of $a_i$
to be the sum over all trees in $0$--dimensional moduli spaces of the
monomial associated to the tree.

We conjecture that the resulting differential graded algebra
$(\A,\d)$ is stable tame isomorphic to the knot DGA from \cite{I},
at least over $\Z_2$; establishing an equivalence over $\Z$ would
entail sorting through orientation issues on the above moduli
spaces, \`a la Morse homology. The knot DGA conjecturally represents
a relative contact homology theory, as described in \cite[Section 3]{I},
which bears a striking resemblance to the DGA described above. In
particular, the DGA of the contact homology is generated by binormal
chords, and the differential is also given by gradient flow trees.
However, the assignment of grading in the two DGAs is different in
general, as is the differential.

We remark that it should be possible to bound the total number of
binormal chords for a knot by examining the full knot contact
homology $HC_*(K)$, similarly to minimal chords and $HC_0$.

\subsection{Cords and $\Sigma_2(K)$}
\label{ssec:sigma2}

In \cite{I}, it was demonstrated that knot contact homology has a
close relation to the double branched cover $\Sigma_2(K)$ of $S^3$
over the knot $K$. Cords can be used to elucidate this relationship,
as we will now see. Recall that the ($SL_2(\C)$) \textit{character
variety} of $\Sigma_2(K)$ is the variety of characters of $SL_2(\C)$
representations of $\pi_1(\Sigma_2(K))$.

\begin{proposition}
There is a map over $\C$ from $HC_0^{\ab}(K)\otimes\C$ to the
coordinate ring of the character variety of $\Sigma_2(K)$.
\label{prop:char}
\end{proposition}

\begin{proof}
Given an $SL_2(\C)$ character $\chi\colon\thinspace
\pi_1(\Sigma_2(K))\to\C$, we wish to produce a map
$HC_0^{\ab}(K)\otimes\C\to\C$. Note that $\chi$ satisfies
$\chi(e)=2$, $\chi(g^{-1})=\chi(g)$, and
$\chi(g_1g_2)+\chi(g_1^{-1}g_2)=\chi(g_1)\chi(g_2)$ for all
$g,g_1,g_2\in\pi_1(\Sigma_2(K))$.

Any (unoriented) cord $\gamma$ of $K$ has two lifts to $\Sigma_2(K)$
with the same endpoints; arranging these lifts head-to-tail gives an
element $\tilde{\gamma}$ of $\pi_1(\Sigma_2(K))$ which is unique up
to conjugation and inversion. In particular, $\chi(\tilde{\gamma})$
is well defined. We claim that the map sending each cord $\gamma$ to
$-\chi(\tilde{\gamma})$ descends to the desired map
$HC_0^{\ab}(K)\otimes\C\to\C$; this simply entails checking the
skein relations in the definition of the cord ring
(Definition~\ref{def:cordring}). The second skein relation
(\ref{eq:skein2}) is preserved since $\chi(e)=2$. As for the first
relation (\ref{eq:skein1}), label the cords depicted in
(\ref{eq:skein1}) by $\gamma_3,\gamma_4,\gamma_1,\gamma_2$ in order,
so that the relation reads $\gamma_3+\gamma_4+\gamma_1\gamma_2=0$.
If we choose the base point for $\pi_1(\Sigma_2(K))$ to be the point
on the knot depicted in the skein relation, then $\tilde{\gamma}_3$
and $\tilde{\gamma}_4$ are conjugate to
$\tilde{\gamma}_1\tilde{\gamma}_2$ and
$\tilde{\gamma}_1^{-1}\tilde{\gamma}_2$ in some order. Since
$\chi(\tilde{\gamma}_1\tilde{\gamma}_2)+
\chi(\tilde{\gamma}_1^{-1}\tilde{\gamma}_2)=\chi(\tilde{\gamma}_1)
\chi(\tilde{\gamma}_2)$, (\ref{eq:skein1}) is preserved.
\end{proof}

If $K$ is two-bridge, then $\Sigma_2(K)$ is a lens space and
$\pi_1(\Sigma_2(K)) \cong \Z/n$ where $n=\Delta_K(-1)$. In this
case, all $SL_2(\C)$ representations of $\pi_1(\Sigma_2(K))$ are
reducible, with character on the generator of $\Z_n$ given by
$\omega^k+\omega^{-k}$ where $\omega$ is a primitive $n$-th root of
unity and $0\leq k\leq n-1$. It follows that the coordinate ring of
the character variety of $\Sigma_2(K)$ is $\C[x]/(p_{(n+1)/2}(x))$
with $p_{(n+1)/2}(x) = \prod_{k=0}^{(n-1)/2} \left(
x-\omega^k-\omega^{-k}\right)$. This is precisely the polynomial
defined inductively in Section~\ref{ssec:2bridge}. It follows from
Theorem~\ref{thm:2bridge} that the map in
Proposition~\ref{prop:char} is an isomorphism when $K$ is
two-bridge.

In fact, a number of calculations on small knots lead us to propose
the following.

\begin{conjecture}
The map in Proposition~\ref{prop:char} is always an isomorphism; the
complexified cord ring $HC_0^{\ab}(K)\otimes\C$ is precisely the
coordinate ring of the character variety of $\Sigma_2(K)$.
\end{conjecture}

In general, the surjection $HC_0^{\ab}(K) \twoheadrightarrow
\Z[x]/(p_{(n(K)+1)/2})$ from \cite[Theorem 7.1]{I}, where $n(K)$ is
the largest invariant factor of $H_1(\Sigma_2(K))$, can be seen via
the approach of Proposition~\ref{prop:char} by restricting to
reducible representations. Proving surjectivity from this viewpoint
takes a bit more work, though.

We can also use cords to see the surjection $HC_0^{\lin,\ab}(K)
\twoheadrightarrow \Sym^2(H_1(\Sigma_2(K)))$ from
\cite[Proposition~7.11]{I}, which was cited in the proof of
Corollary~\ref{cor:minchords} above. Consider $H_1(\Sigma_2(K))$ as
a $\Z$--module, with group multiplication given by addition. Given a
cord of $K$, we obtain an element of $H_1(\Sigma_2(K))$, as in the
proof of Proposition~\ref{prop:char}, defined up to multiplication
by $\pm 1$; the square of this element gives a well-defined element
of $\Sym^2(H_1(\Sigma_2(K)))$. The fact that this map descends to
$HC_0^{\lin,\ab}(K)$ is the identity $(x+y)^2+(x-y)^2-2x^2-2y^2=0$.
Again, proving that this map is surjective takes slightly more work.

\subsection{Extensions of the cord ring}
\label{ssec:extensions}

We briefly mention here a couple of extensions of the knot cord
ring. These each produce new invariants which may be of interest.

One possible extension is to define a cord ring for any knot in any
$3$--manifold, in precisely the same way as in $\R^3$. If the
$3$--manifold and the knot are sufficiently well-behaved (eg, no
wild knots), it seems likely that the cord ring will always be
finitely generated. The cord ring would be a natural candidate for
the degree $0$ portion of the appropriate relative contact homology
\cite[Section 3]{I}. It would be interesting to construct tools to
compute the cord ring in general, akin to the methods used in
\cite{I} and this paper. See also the Appendix.

As another extension, we note that the definition of the cord ring
makes sense not only for knots, but also for graphs embedded in
$\R^3$, including singular knots. It is clear for topological
reasons that the cord ring is invariant under neighborhood
equivalence; recall that two embedded graphs are neighborhood
equivalent if small tubular neighborhoods of each are ambient
isotopic.

\begin{proposition}
\label{prop:handlebody}
The cord ring is an invariant of graphs embedded in $\R^3$, modulo
neighborhood equivalence.
\end{proposition}

By direct calculation, one can show that the cord ring is a
nontrivial invariant for graphs of higher genus than knots. For
instance, the graph consisting of the union of a split link and a
path connecting its components has cord ring which surjects onto the
cord ring of each component of the link. By contrast, the figure
eight graph (or the theta graph), like the unknot, has trivial cord
ring $\Z$.

The graph cord ring can be applied to tunnel numbers of knots, via
the observation that any graph whose complement is a handlebody has
trivial cord ring. This can be used to compute lower bounds for the
tunnel numbers of some knots, but the process is somewhat laborious.


\section*{Appendix: The cord ring and fundamental groups}
\addcontentsline{toc}{section}{Appendix: The cord ring and fundamental groups}

\renewcommand{\thesection}{A}
\setcounter{proposition}{0}

In this appendix we show that the cord ring is determined by the
fundamental group and peripheral structure of a knot. We then
introduce a generalization of the cord ring to any codimension $2$
submanifold of any manifold and derive a homotopy-theoretic
formulation in this more general case. As an application, we show
that the cord ring gives a nontrivial invariant for embeddings of
$S^2$ in $S^4$.

Assume that we are given a knot $K\subset S^3$; it is
straightforward to modify the constructions here to links. Let
$N(K)$ denote a regular neighborhood of $K$ and let $M=S^3\setminus
\int(N(K))$ denote the knot exterior. Write $G=\pi_1(M)$, and let $P$
be the image of $\pi_1(\d N(K))$ in $G$; define $\CC(G,P) =
P\backslash G/P$, and write $[g]$ for the image of $g\in G$ in
$\CC(G,P)$.

\begin{lemma}
The set $\CC_K$ of homotopy classes of cords of $K$ is identical to
$\CC(G,P)$.
\label{lem:coset}
\end{lemma}

\begin{proof}
As in the proof of the Van Kampen theorem, it is easy to see that
there is a one-to-one correspondence between homotopy classes of
cords for $K$ and homotopy classes of cords in the knot exterior
$M$, where a cord in $M$ is a continuous path
$\alpha\colon\thinspace [0,1]\to M$ with $\alpha^{-1}(\d
M)=\{0,1\}$. We will identify the latter set of homotopy classes
with $P\backslash G/P$.

Fix a base point $x_0$ in $\d M$. Given a cord $\alpha$ in $M$, we
pick paths $\beta$ and $\gamma$ in $\d M$ joining $x_0$ to
$\alpha(0)$ and $\alpha(1)$. We associate to $\alpha$ the
equivalence class $[\beta\alpha\gamma^{-1}]\in P\backslash G/P$.
This is clearly independent of the choice of $\beta$ and $\gamma$.
Furthermore, homotopic cords give the same element of $P\backslash
G/P$, because a family of cords can be given a continuous family of
paths $\beta$, $\gamma$. Hence we have a map from the set of cords
of $M$ to $P\backslash G/P$.

To construct the inverse of this map, observe that each element in
$G$ has a representative which does not intersect $\d M$ in its
interior, and hence gives a cord which is unique up to homotopy; in
addition, for any $g\in G$ and $h,k\in P$, $g$ and $hgk$ give
homotopic cords. This completes the proof of the lemma.
\end{proof}

Using Lemma~\ref{lem:coset}, we can reformulate the definition of
the cord ring in group-theoretic terms. Let $\A(G,P)$ be the tensor
algebra freely generated by the set $\CC(G,P)$, let $\mu\in G$
denote the homotopy class of the meridian of $K$, and let $I(\mu)$
be the ideal in $\A(G,P)$ generated by the ``skein relations''
\[
[\alpha\mu\beta]+[\alpha\beta]+[\alpha]\cdot [\beta],~~~~~
\alpha,\beta\in G,
\]
and $[e]+2$, where $e=1\in G$.

\begin{proposition}
$\A(G,P)/I(\mu)$ is the cord ring of the knot $K$.
\end{proposition}

\begin{proof}
The skein relations generating $I(\mu)$ are simply the
homotopy-theoretic versions of the skein relations in the cord ring.
\end{proof}

Note that the above construction associates a ring to any triple
$(G,P,\mu)$, where $G$ is a group, $P$ a subgroup, and $\mu$ an
element of $P$. Such a triple is naturally associated to any
codimension $2$ embedding $K\subset M$ of manifolds; we will be more
precise presently. In this general setting, we can introduce a cord
ring which agrees with $\A(G,P)/I(\mu)$, and which specializes to
the usual cord ring for knots in $S^3$.

\begin{definition}
Let $K\subset M$ be a codimension $2$ submanifold. A \textit{cord}
of $K$ is a continuous path $\gamma\colon\thinspace [0,1]\rightarrow
M$ with $\gamma^{-1}(K) = \{0,1\}$. A \textit{near homotopy of
cords} is a continuous map $\eta\colon\thinspace [0,1]\times
[0,1]\rightarrow M$ with $\eta^{-1}(K) = ([0,1] \times \{0\}) \cup
([0,1] \times \{1\}) \cup \{(t_0,s_0)\}$, for some $(t_0,s_0) \in
(0,1) \times (0,1)$ such that $\eta$ is transverse to $K$ in a
neighborhood of $(t_0,s_0)$.
\end{definition}

\noindent Less formally, a near homotopy of cords is a homotopy of
cords, except for one point in the homotopy where the cord breaks
into two.

Just as for knots in $S^3$, let $\CC_K$ denote the set of homotopy
classes of cords of $K$, and let $\A_K$ be the tensor algebra freely
generated by $\CC_K$. To each near homotopy of cords, we can
associate an element in $\A_K$, namely $[\gamma_0] + [\gamma_1] +
[\gamma_2] \cdot [\gamma_3]$, where
$\gamma_0,\gamma_1,\gamma_2,\gamma_3$ are the cords of $K$
corresponding to $\eta(\{0\}\times [0,1])$, $\eta(\{1\}\times
[0,1])$, $\eta(\{t_0\} \times [0,s_0])$, $\eta(\{t_0\} \times
[s_0,1])$, respectively. Now define $\I_K$ to be the ideal in $\A_K$
generated by the elements associated to all possible near
homotopies, along with the element $[e]+2$, where $e$ represents the
homotopy class of a contractible cord. (For a knot in $\R^3$, this
agrees with the skein-relation definition of $\I_K$ used to
formulate the original cord ring.)

\begin{definition}
The \textit{cord ring} of $K \subset M$ is $\A_K/\I_K$.
\end{definition}

It is clear that the cord ring is an invariant under isotopy. We
have seen that, for knots in $S^3$, the cord ring can be written
group-theoretically, in terms of the peripheral structure of the
knot group. A similar expression can be given for the cord ring of a
general codimension $2$ submanifold $K \subset M$. Let $N(K)$ denote
a tubular neighborhood of $K$ in $M$; its boundary is a circle
bundle over $K$. Set $G = \pi_1(M\setminus K)$ with base point $p$
on $\d N(K)$, $P = i_*\pi_1(\d N(K))$ where $i$ is the inclusion $\d
N(K) \hookrightarrow M\setminus K$, and $\mu$ equals the homotopy
class of the $S^1$ fiber of $\d N(K)$ containing $p$.

\begin{proposition}
For any codimension $2$ submanifold $K$, the cord ring of $K$ is
isomorphic to $\A(G,P)/I(\mu)$.
\end{proposition}

\begin{proof}
Completely analogous to the proof for knots in $S^3$.
\end{proof}

We now consider a particular example of the cord ring, for
embeddings of $S^2$ in $S^4$. Recall that any knot in $S^3$ yields a
``spun knot'' $2$--sphere in $S^4$; see, eg, \cite{Rol}.

\begin{proposition}
The cord ring distinguishes between the unknotted $S^2$ in $S^4$ and
the spun knot obtained from any knot in $S^3$ with nontrivial cord
ring (in particular, any knot with determinant not equal to $1$).
\end{proposition}

\begin{proof}
In the case of the unknotted $S^2$, $G=P=\Z$ and hence the cord ring
is trivial. On the other hand, suppose that $K$ is a knot with
nontrivial cord ring. For the spun knot obtained from $K$, $G$ is
$\pi_1(S^3\setminus K)$, $\mu$ is the element corresponding to the
meridian of $K$, and $P$ is the subgroup of $G$ generated by $\mu$.
It follows that the cord ring of the spun knot surjects onto the
cord ring for $K$, and hence is nontrivial.
\end{proof}

\noindent Thus the cord ring gives a nontrivial invariant for a
large class of $2$--knots in $S^4$.

Just as the cord ring for knots in $\R^3$ should give the
zero-dimensional relative contact homology of a certain Legendrian
torus in $ST^*\!\R^3$, we believe that the cord ring in general
should correspond to a zero-dimensional contact homology. Recall
from, eg, \cite{I} that any submanifold $K\subset M$ gives a
Legendrian submanifold $LK$ of the contact manifold $ST^*\!M$ given
by the unit conormal bundle to $K$.

\begin{conjecture}
For any codimension $2$ submanifold $K \subset M$, the cord ring of
$K$ is the zero-dimensional relative contact homology of $LK$ in
$ST^*\!M$.
\end{conjecture}

Another natural direction of inquiry is to consider
higher-dimensional contact homology for knots. Viterbo (\cite{Vit},
see also \cite{AS,SW}) has shown that the Floer homology of the
tangent bundle of a manifold is the cohomology of its loop space.
Here we have shown how the zero-dimensional contact homology of a
knot can similarly be determined in terms of the algebraic topology
of the space of cords. It seems possible that higher-dimensional
contact homology may have a description analogous to our
description, except that one takes into account not just the
homotopy classes of cords, but the full homotopy type of the space
of cords.



\begin{thebibliography}{10}

\bibitem{AS}
\textbf{A Abbondandolo}, \textbf{M Schwarz}, \textit{On the Floer
homology of cotangent bundles}, Comm.\ Pure Appl.\ Math. to appear

\bibitem{Big}
\textbf{S\,J Bigelow}, \emph{Braid groups are linear}, J. Amer. Math. Soc. 14
  (2001) 471--486 \MR{1815219}

\bibitem{Bir}
\textbf{J\,S Birman}, \emph{Braids, links, and mapping class groups}, Princeton
  University Press, Princeton, N.J. (1974) \MR{0375281}

\bibitem{Hum}
\textbf{S\,P Humphries}, \emph{An approach to automorphisms of free groups and
  braids via transvections}, Math. Z. 209 (1992) 131--152 \MR{1143219}

\bibitem{Kra}
\textbf{D Krammer}, \emph{The braid group {$B\sb 4$} is linear}, Invent. Math.
  142 (2000) 451--486 \MR{1804157}

\bibitem{Mag}
\textbf{W Magnus}, \emph{Rings of {F}ricke characters and automorphism groups
  of free groups}, Math. Z. 170 (1980) 91--103 \MR{558891}

\bibitem{I}
\textbf{L Ng}, \emph{Knot and braid invariants from contact homology. {I}},
  \gtref9{2005}8{247}{297} \MR{MR2116316}

\bibitem{Rol}
\textbf{D Rolfsen}, \emph{Knots and links}, Publish or Perish Inc., Berkeley,
  Calif. (1976) \MR{0515288}

\bibitem{SW}
\textbf{D Salamon}, \textbf{J Weber}, \textit{Floer homology and the
heat flow}, e-print \arxiv{math.SG/0304383}

\bibitem{Vit}
\textbf{C Viterbo}, \textit{Functors and computations in Floer
homology with applications II}, preprint 1998

\end{thebibliography}
\end{document}